\documentclass{article}

\usepackage{amsmath,amssymb}

\usepackage{epsf}
\usepackage{amsfonts}
\usepackage{amstext}
\usepackage{amsbsy}
\usepackage{amsthm}
\usepackage{amscd}
\usepackage{graphics}
\usepackage{rotating}
\usepackage[all]{xy}

\newtheorem{Thm}{Theorem}
\newtheorem{Lem}{Lemma}
\newtheorem{Cor}{Corollary}

\begin{document}
\title{Noncyclic Covers of Knot Complements}
\author{
   Nathan Broaddus \footnote{The author is supported by
an NSF Postdoctoral Fellowship at Cornell University} \\
   Cornell University \\
}
\date{March 22, 2005}
\maketitle

\begin{abstract}
Hempel has shown that the fundamental groups of knot complements
are residually finite. This implies that every nontrivial knot must
have a finite-sheeted, noncyclic cover. We give an explicit bound,
$\Phi (c)$, such that if $K$ is a nontrivial knot in the three-sphere
with a diagram with $c$ crossings and a complement with a particularly
simple JSJ decomposition then the complement of $K$ has a
finite-sheeted, noncyclic cover with at most $\Phi (c)$ sheets. 
\end{abstract}






\section{Introduction}
\label{introduction_section}

Let $K$ be a nontrivial knot in $S^3$.  Let $M = S^3 - N(K)$ 
be the complement of an open regular neighborhood of $K$.
It is well known that for each positive integer, 
$k$, $M$ has a unique cyclic cover with $k$ sheets arising from the map 
$\pi_1(M) \to H_1(M) \cong \mathbb{Z} \to \mathbb{Z}/k\mathbb{Z}$.
However, much less is known 
about the noncyclic covers of knot complements.  In \cite{Hempel} 
Hempel establishes that fundamental groups of Haken 3-manifolds are 
residually finite. This shows in particular that the fundamental 
groups of knot complements are residually finite.  Thus for any 
nontrivial element, $g$, of the commutator of $\pi_1(M)$ there must be a 
nontrivial normal subgroup of finite index in $\pi_1(M)$ not containing 
$g$. It follows that knot complements must have infinitely many 
finite, nonabelian covers. The goal of this exposition is to give
an explicit function $\Phi (c)$ such that
if $K$ is a nontrivial knot with a diagram with $c$ crossings
and its complement $M = S^3 - N(K)$ has a simple JSJ decomposition 
then $M$ has a noncyclic cover with at most $\Phi (c)$ 
sheets.

The question behind this investigation is how well finite index subgroups
differentiate the fundamental group of a nontrivial knot complement from
the group of integers.  The bound on the index given here
(see Theorem \ref{mainthm})
seems to be largely an artifact of the techniques used.  It is safe
to conjecture that it might be drastically improved with another approach.

If one could generalize these results to all knots then in a technical
sense, one would get an algorithm
for detecting knottedness as follows:  If one starts with
a knot with $c$ crossings and systematically creates all 
covers of the complement with $\Phi(c)$ or less sheets then if a 
noncyclic cover is found, the knot is nontrivial.  If no such cover is 
found, the knot is trivial.
However, in light of the large bound given in this paper, much better,
if still impractical algorithms already exist to establish knottedness.
In \cite{Hass}, a bound on the number of Reidemeister moves needed to
convert an arbitrary diagram of the unknot to the standard diagram of the
unknot is given.  Also, Ian Agol has shown that computing lower
bounds for the genus of a knot is in NP. At this writing, it is an open
problem whether there is a practical algorithm to detect knottedness.

\vspace{0.2in}

\noindent {\bf \large Acknowledgments}

\vspace{0.15in}

This paper is based on the author's Ph.D. thesis at Columbia University.
The author would like to thank his advisor, Joan Birman, for her guidance,
support, and encouragement.  The author would also like to thank Columbia
University and the National Science Foundation Graduate Fellowship Program
for support of his studies.  The author owes a debt of thanks to Walter
Neumann, David Bayer, Shou-Wu Zhang, Ian Agol, Saul Schleimer, Brian Mangum,
Douglas Zare, Tara Brendle, and Abhijit Champanekar for helpful
suggestions and conversations.  Special thanks are extended to
Xingru Zhang, Ryan Budney and the anonymous referee for
pointing out errors in a previous version.




\section{Main Result}
\label{main_result_section}

Set
\begin{equation}
A(n) =\frac{(n^2-n+1)!}{n^2[(n-1)!]^n},
\label{adef}
\end{equation}
\begin{equation}
 B(n)=\frac{n^3-n^2}{n^2-n+1},
 \label{bdef}
\end{equation}
\begin{eqnarray}
\lefteqn{D(n) =} \nonumber \\
& & \exp
        \Big[ 2(4n+4)\left( 8n^2 + 4n \right)^{2^{4n+4}}
A(4n+5)( 27n + 5) \nonumber \\
& & \cdot \left( 2^{4n+2} + 3\cdot 2^{3n+3} + \left({\scriptstyle 
\frac{\sqrt{3}}{2}} + 3 \right) n 
+B(4n+5)\log 2  +{\scriptstyle \frac{\sqrt{3}}{2}} \right) \nonumber \\
& & + 2^{4n+4} \left( 8n^2 + 4n \right)^{(4n+4)2^{4n+4}}
         \Big( 2\log 2 + 4 \log (16n^2 3^{n-1}) \nonumber \\ [-3mm]
& & \label{primebound} \\ [-3mm]
& & + 3(2^{4n-1}-1)( \log 2 ) + 3(2^{4n}+2^{4n-1}-2)A(4n+5)( 27n + 5) 
\nonumber \\
& & \cdot \left( 2^{4n+2} + 3\cdot 2^{3n+3} + \left({\scriptstyle 
\frac{\sqrt{3}}{2}} + 3 \right) n 
+B(4n+5)\log 2  +{\scriptstyle \frac{\sqrt{3}}{2}} \right) \Big) \Big]
        \nonumber
\end{eqnarray}
and
\begin{equation}
\Phi (c) = \left(87 \left(\log (D(100c)) + 8c \log {\textstyle
\frac{c}{2}} \right) \right)^{24c \left( 2^{4n+4} \left(8n^2 + 4n \right)^{
(4n+4)2^{4n+4}}\right)}.
\end{equation}
A knot $K$ in $S^3$ will be called {\em decompositionally linear} if
the JSJ decomposition of its complement along essential tori $T_1, T_2,\cdots, T_r$
has the property that in $S^3$, $K$ and $T_i$ are on the same side of $T_j$ if
$j < i$.

The main result is as follows:
\begin{Thm}
Let $K$ be a nontrivial, decompositionally linear knot in $S^3$ and $M = S^3 - K$ its complement. 
Suppose $K$ has a diagram with $c$ crossings. Then $M$ has a noncyclic
cover with at most $\Phi (c)$ sheets.
\label{mainthm}
\end{Thm}

The proof of Theorem \ref{mainthm} will proceed using Thurston's 
geometrization for knot complements. Let $K$ and $M$ be as in the 
statement of the theorem. Then the JSJ decomposition of $M$ cuts $M$ 
along essential tori $T_1, T_2,\cdots, T_r$ into spaces $M_0, M_1, 
\cdots, M_r$ where either $M_i$ is Seifert fibered or $M_i - \partial 
M_i$ has a complete hyperbolic structure.




\section{Topology of knot complements}
\label{topology_of_knot_complements_section}




\subsection{Standard spines and ideal triangulations}
\label{standard_spines_and_ideal_triangulations_subsection}

For our purposes an \emph{ideal triangulation} of a $3$-manifold, $M$, will
be a simplicial complex, $\mathcal{T}$, satisfying some further
conditions.  The complex, $\mathcal{T}$, must be a union of a finite number
of $3$-simplices with pairs of faces identified.  In fact, we
insist that there are no unidentified ``free'' faces.
Identification of different faces of the same tetrahedron
will be allowed. For $\mathcal{T}$ to be an ideal triangulation of
the $3$-manifold, $M$,
we require $\mathcal{T}$ minus its vertices to be homeomorphic to
$M-\partial M$.  We will write
$\mathcal{T} = \bigcup_{i=1}^n \sigma_i$ to indicate that
$\mathcal{T}$ is an ideal triangulation with the $n$ ideal tetrahedra,
$\sigma_1, \sigma_2, \dots, \sigma_n$.
For our purposes, the links of the vertices in our ideal
triangulations will always be tori, and all $3$-manifolds and
their ideal triangulations will be orientable.

When dealing with the geometric pieces of $M$ it will
be convenient to have a bound on the number of tetrahedra
needed to triangulate them.  In working with ideal triangulations
it is often helpful to be familiar with the dual notion of
standard spines. As in \cite{Casler} a \emph{spine} is simply
a $2$-complex. The \emph{singular 1-skeleton}
of a spine is the set of points which do not have neighborhoods
homeomorphic to open disks. The \emph{singular vertices}
of a spine are the points of the singular 1-skeleton
which do not have neighborhoods in the singular 1-skeleton
homeomorphic to open intervals.

Let $C$ be a spine, $C_1$ be the singular 1-skeleton 
of $C$, and $C_0$ be the singular vertices of $C$. The spine
$C$ will be a \emph{standard spine} if it satisfies three
conditions.  Firstly, $C$ must satisfy the \emph{neighborhood
condition}.  That is, every point of $C$ must have a neighborhood
homeomorphic to one of the three 2-complexes
\begin{figure}[ht]
$$
\setlength{\unitlength}{0.05in}
\includegraphics{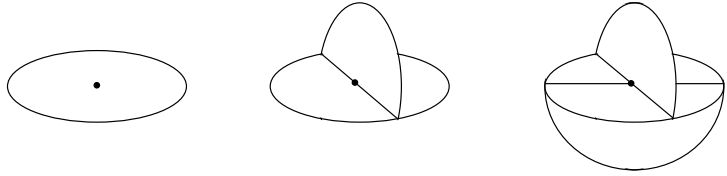}
$$
\caption{The three possible neighborhoods in a standard spine}
\label{spine_fig}
\end{figure}
pictured in Figure \ref{spine_fig}. Secondly, $C - C_1$ must be a union
of countably many disjoint 2-disks.  Thirdly, we require that
$C_1 - C_0$ be a union of countably many disjoint arcs.
The complex, $C$, is a \emph{spine} (resp. \emph{standard spine})
\emph{of a
3-manifold} $N$ if $C \subset N$ is a spine (resp. standard spine) and 
$N$ collapses to $C$.
An important property of a standard spines is that
if $C$ is a standard spine of $N$ then if $C$ is embedded
in any $3$-manifold then $N$ is homeomorphic
to a regular neighborhood of $C$ in that manifold.

\begin{figure}[ht]
$$
\setlength{\unitlength}{0.05in}
\includegraphics{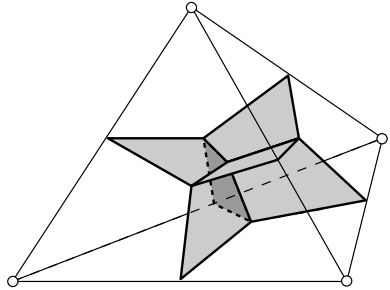}
$$
\caption{An ideal tetrahedron and its dual spine}
\label{dual_fig}
\end{figure}

In \cite{TV} and \cite{Petronio} it is mentioned that standard
spines are dual to ideal triangulations (see Figure \ref{dual_fig}).
For every ideal triangulation of a 3-manifold there is a dual
standard spine, and for every standard spine there is a dual
ideal triangulation.  Thus we see that a standard spine
carries the same information as an ideal triangulation.
Moreover, the number of singular
vertices in a standard spine will be the number of
ideal tetrahedra in the dual triangulation.  We will exploit this
duality a number of times.




\subsection{Triangulating a knot complement}
\label{triangulating_a_knot_complement_subsection}

A preliminary step in our exposition will be to relate
the number of crossings in a knot diagram to the number
of ideal tetrahedra needed to triangulate the complement
of the knot.  It is noted in \cite{Petronio} that
the number of ideal tetrahedra needed is at most linear
in the number of crossings in a projection.  Here we
give an explicit relationship.  The argument is essentially based
on the triangulation algorithm in Jeff Weeks' program,
SnapPea.
\begin{Lem}
Let $K$ be a knot in $S^3$ with a diagram with $c > 0$ crossings.  Then
the complement of $K$ has an ideal triangulation with less than
$4c$ ideal tetrahedra.
\label{complement_triangulation_lemma}
\end{Lem}
A proof of this lemma is given in Appendix
\ref{complement_triangulation_appendix}.

Now we translate the bound on the number of ideal tetrahedra
needed to triangulate $M$ into a bound on the
number needed to triangulate the geometric pieces of $M$.

\begin{Lem}
Let $K$ be a knot in $S^3$ and $M = S^3 - K$ its complement.  
Suppose $M$ can be triangulated with $t$ ideal tetrahedra.  
Also suppose that embedded, disjoint tori, $T_1,T_2, \cdots, T_r$, 
give the JSJ decomposition of $M$, and $M_0, M_1, \cdots, M_r$ are the 
connected components after cutting. Then the $M_i$'s have ideal
triangulations with $t_i$ ideal tetrahedra each so that $\sum_{i=0}^r 
t_i \leq 25t$.
\label{decomposition_triangulation_lemma}
\end{Lem}

\begin{proof}
Let $K$ and $M$ be as in the statement of the lemma.  By assumption
$M$ has an ideal triangulation $\mathcal{T} = \bigcup_{i=1}^t
\sigma_i$ with $t$ ideal tetrahedra.  We may choose
our tori $T_1, T_2, \cdots, T_r$ so that their union $S$ is a normal 
surface with respect to $\mathcal{T}$.  For each $i, 1 \leq i \leq t$, 
$S$ cuts
$\sigma_i$ into pieces with four basic types (see Figure 
\ref{type1234_fig}):
\begin{itemize}
\item[{\bf (a)}] Pieces whose closure intersects $S$ in 
two triangles. \label{chunk1}
\item[{\bf (b)}] Pieces whose closure intersects $S$ in two 
quadrilaterals. \label{chunk2}
\item[{\bf (c)}] Pieces whose closure intersects $S$ in 
two triangles and one quadrilateral. \label{chunk3}
\item[{\bf (d)}] Pieces whose closure intersects $S$ in 
four triangles. \label{chunk4}
\end{itemize}
Of course there will also be pieces whose closures in $\sigma_i$
will be incident with 
the corners of $\sigma_i$, but we will put these pieces in categories
{\bf (a)} - {\bf (d)} based on how they look
when the corners of $\sigma_i$ are cut off by triangles.

We will now construct a standard spine of $M-S$. Consider the $i$th
tetrahedron in the ideal triangulation of $M$ and its intersection
with the surface $S$.  For each region of type
{\bf (a)} and {\bf (b)} place a triangular or quadrilateral disk in 
its center parallel to the faces incident with $S$ as 
shown in Figure \ref{type1234_fig}.
\begin{figure}[ht]
$$
\setlength{\unitlength}{0.05in}
\put(9.0,41.5){\footnotesize Type {\bf (a)}}
\put(38.0,41.5){\footnotesize Type {\bf (b)}}
\put(11.4,3.7){\footnotesize $S$}
\put(11.6,0.5){\footnotesize spine}
\put(9.0,10){\footnotesize Type {\bf (c)}}
\put(38.0,10){\footnotesize Type {\bf (d)}}
\includegraphics[scale=0.7]{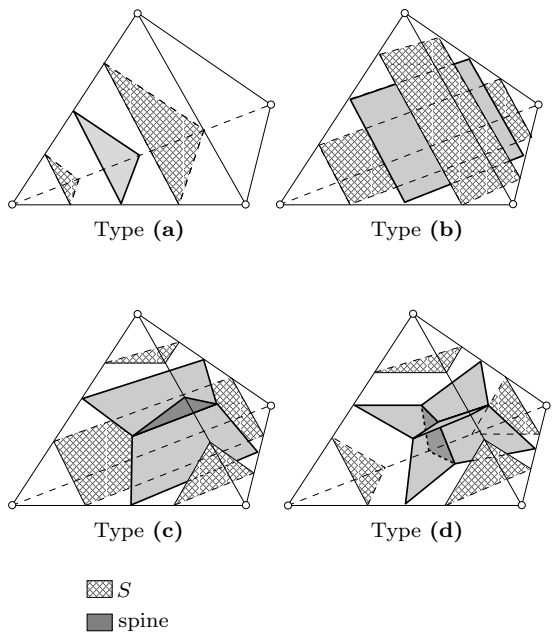}
$$
\caption{Regions of type {\bf (a)} - {\bf (d)} and their spines}
\label{type1234_fig}
\end{figure}
For each region of type {\bf (c)} and {\bf (d)} place a 2-complex as 
in Figure \ref{type1234_fig}.
Let $D_i \subset \sigma_i$ be the union of all of these 
spine pieces.  One sees immediately that 
$C' = \bigcup_{i = 1}^t D_i$ is a spine of $M-S$ and
that it satisfies the neighborhood condition.  The spine, $C'$,
has a singular vertex for each region of type {\bf (d)}. Clearly 
each tetrahedron of $\mathcal{T}$ contains at most one region of type 
{\bf (d)}. Thus $C'$
has at most $t$ vertices.  Note that each connected component of $C'$ 
has a nonempty singular $1$-skeleton, for if there were a component, $A$, of 
$C'$ with empty singular $1$-skeleton then the component of $M-S$
containing $A$ would be homeomorphic to the torus crossed with the
open unit interval.

Let $C_1'$ be the singular 1-skeleton of $C'$.  Although $C'$ satisfies
the neighborhood condition, $C'$ will not in general be a standard
spine since $C' - C_1'$ may not be a union of disks.   A component of $C' 
- C_1'$ must have genus 1 or 0 since the boundary components of $M-S$ 
are tori.  Thus if a component of $C' - C_1'$ has $b$ boundary 
components it may be cut into disks with $b+1$ or less arcs. For each 
of these arcs, $\gamma$, modify $C'$ as in Figure \ref{arcmod_fig} (This is 
possible since each connected component of $C'$ has a nonempty singular 
$1$-skeleton).
\begin{figure}[ht]
$$
\setlength{\unitlength}{0.05in}
\put(9,10){\footnotesize \bf $C'$}
\put(46,10){\footnotesize \bf $C$}
\includegraphics{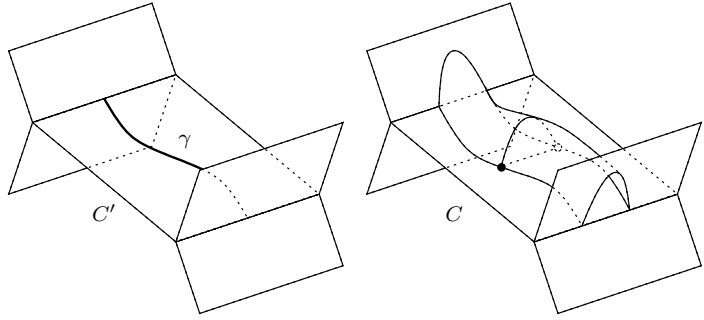}
\put(-55,18){\footnotesize \bf $\gamma$}
$$
\caption{Modifying $C'$ to get $C$}
\label{arcmod_fig}
\end{figure}
Let $C$ be the modified spine.  The modification in Figure 
\ref{arcmod_fig} takes one spine to
another \cite{Piergallini} so $C$ is a spine of $M-S$.  The spine, $C$,
satisfies the neighborhood condition. Let $C_1$ be the singular
1-skeleton of $C$ and $C_0$ be the
singular vertices of $C_1$. By construction, $C - C_1$ is a disjoint
union of disks.  We claim that $C_1 - C_0$ must be a disjoint union of
arcs. If not then $C_1 - C_0$ contains an $S^1$. Since $C - C_1$ is composed
entirely of disks, each disk in the component of $C$ containing this
$S^1$ must have this $S^1$ as its boundary.
This is impossible because, as mentioned before, the only such spine
satisfying the neighborhood
condition is three disks glued along their boundary. This is not a spine
of the complement of a knot in the $3$-sphere or the complement
of a knot in a solid torus. Consequently, $C$ is a
standard spine of $M-S$.

Let us now count the singular vertices of $C$. The total number of 
boundary components of $C'-C'_1$ is less than or equal to 3 times the 
number of regions of type {\bf (c)} plus 6 times the number of regions 
of type {\bf (d)}.  The number of components of $C'-C'_1$ is bounded by 
that same number.  Let $\mathcal{A}$ be the set of connected
components of $C'-C'_1$, and for $A \in \mathcal{A}$ let $b_A$ be
the number of boundary components of $A$. If $a$ is the number of
arcs needed to cut 
$C' - C'_1$ into disks we have:
\begin{eqnarray*}
a & = & \sum_{{A \in \mathcal{A}} \atop {\mathrm{g}(A) = 1}} (b_A +1)
+ \sum_{{A \in \mathcal{A}} \atop {\mathrm{g}(A) = 0}} (b_A - 1) \\
& = & \sum_{A \in \mathcal{A}} b_A + \sum_{{A \in \mathcal{A}}
\atop {\mathrm{g}(A) = 1}} 1
- \sum_{{A \in \mathcal{A}} \atop {\mathrm{g}(A) = 0}} 1 \\
& \leq & 6t + | \mathcal{A} | \\
& \leq & 12t
\end{eqnarray*}
Changing $C'$ to $C$ introduces 2 singular vertices for each cutting 
arc so $C$ has at most $2a \leq 2 \cdot 12t = 24t$ more singular vertices than 
$C'$.  As mentioned above, $C'$ has at most $t$ vertices so $C$ has at 
most $25t$ vertices.

The standard spine, $C$, is dual to an ideal
triangulation of $M-S$ with the same number of 
ideal tetrahedra as singular vertices of $C$.  This shows that $M-S$ 
can be triangulated with $25t$ or less tetrahedra. 
\end{proof}

In proof of Lemma \ref{decomposition_triangulation_lemma} we saw that for each $i$,
$M_i \cap C'$ must have nonempty singular 1-skeleton.  The singular
1-skeleton of $C'$ can have, at most, 2 components for every tetrahedron
in $\mathcal{T}$ so we get the following corollary:

\begin{Cor}
Let $K$ be a knot in $S^3$ and $M = S^3 - K$ its complement.  
Suppose $M$ can be triangulated with $t$ ideal tetrahedra.  
Also suppose that the embedded, disjoint tori $T_1,T_2, \cdots, T_r$ 
give the JSJ decomposition of $M$.  Then $r < 2t$.
\label{deccor}
\end{Cor}




\subsection{The JSJ decomposition of certain knot complements}
\label{decomposition_of_a_knot_complement_subsection}

\begin{figure}[ht]
$$
\setlength{\unitlength}{0.05in}
\put(0,3){\footnotesize \bf $M_0$}
\put(10,3){\footnotesize \bf $M_1$}
\put(20,3){\footnotesize \bf $M_2$}
\put(30,3){\footnotesize \bf $M_3$}
\put(45,3){\footnotesize \bf $M_{r-1}$}
\put(56,3){\footnotesize \bf $M_r$}
\put(38,1){\footnotesize \bf $\cdots$}
\put(5.5,-2){\footnotesize \bf $T_1$}
\put(15.5,-2){\footnotesize \bf $T_2$}
\put(25.5,-2){\footnotesize \bf $T_3$}
\put(51.5,-2){\footnotesize \bf $T_r$}
\includegraphics{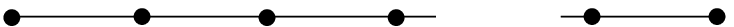}
$$
\caption{$M$ as a graph product}
\label{graph_fig}
\end{figure}

In general the JSJ decomposition of a knot complement can
be quite complicated.  For this study we will restrict to a class
of knots whose complements have a particularly nice JSJ decomposition.
As before, let $K$ be a knot and $M = S^3 - K$ its complement.
Let $T_1,T_2, \cdots, T_r \subset M \subset S^3$ be the tori in the decomposition
for $M$. We restrict $K$ (and possibly re-index the tori) so that
$K$ and $T_i$ are on the same side of $T_j$ whenever $j < i$.
These knots will be called {\em decompositionally linear.}  
Thus we get $M$ as a graph product 
of CW complexes based on the graph in Figure \ref{graph_fig} (See
\cite{Hempel} for more on graph products).

In order to create a covering space of $M$, we will produce a 
compatible collection $\{ \widetilde{M}_i \}$ of finite covers 
of each of the $M_i$'s and assemble them into a finite cover,
$\widetilde{M}$, of $M$ following \cite{Hempel} (see Figures \ref{cover_fig}
and \ref{assemble_fig}).  More specifically, we will choose a prime,
$p$, and let $P \unlhd \mathbb{Z} \times \mathbb{Z}$ be the
characteristic subgroup generated by $(p,0)$ and $(0,p)$.
For each torus, $T_k,$ we will let $\widetilde{T}_k$
be the cover associated to the subgroup of
$\pi_1 (T_k) \cong \mathbb{Z} \times \mathbb{Z}$ corresponding to $P$.
We will then produce a finite cover, $\widetilde{M}_i$, of each
$M_i$ all of whose boundary components will be equivalent to the
appropriate $\widetilde{T}_k$.  Finally we
will assemble copies of these $\widetilde{M}_i$'s to get a cover of $M$.
For each $M_i$ the challenges will be to discover for which
primes, $p$, we will be able to produce such a cover and then
to bound the number of sheets in that cover.

At this point we fix some notation which we use for the rest of
the discussion.  As mentioned above $M_i$ has one or two boundary
components which we denote $\partial_0 M_i = T_{i+1}$ and possibly
$\partial_1 M_i = T_i$.  Recall that $M_i$ is a subset of $S^3$.
For $k \in \{0,1\}$ let $E_k( M_i)$ be the closure of the component of
$S^3 - \partial_k M_i$ disjoint from $M_i$.  
Note that $E_0(M_i)$ is a solid torus,
and $E_1(M_i)$ is a knot complement if it is defined.
A \emph{meridian}, $m_k$, of $M_i$ is an essential,
simple, closed curve in $\partial_k M_i$ with
$m_k$ homologically trivial in $E_k(M_i).$
A \emph{longitude}, $l_k$, is an essential, simple, closed
curve in $\partial_k M_i$ intersecting
a meridian once with the added property that a parallel copy of
$l_k$ in the interior of $M_i$ has linking number 0
with $l_k$ in $S^3$.  When convenient, we will assume that $m_k$
and $l_k$ are oriented loops.




\section{Covers from homology}
\label{covers_from_homology_section}

Let $N = M_i$ be some piece in the JSJ decomposition of the knot
complement, $M$.
In many cases we get an appropriate finite sheeted
covering space $\widetilde{N}$ for $N$ from homology.
The following lemma addresses these cases.

\begin{Lem}
Suppose $N$ is the complement of an open, regular neighborhood of
a knot $L$ in a solid torus.  Suppose further that $L$ has winding
number $w \neq 0$ in the solid torus.  Then for every prime $p$ not dividing
$w$, $N$ has a $p^2$-sheeted covering space $\widetilde{N}$, such that
each boundary component of $\widetilde{N}$ is the cover of a boundary
component of $N$ corresponding to the subgroup $\langle (p,0),(0,p)
\rangle < \mathbb{Z} \times \mathbb{Z} \cong \pi_1 (\partial N ).$
\label{homology_lemma}
\end{Lem}

\begin{proof}
As in the statement of the lemma, let $N$ be the complement of
an open, regular neighborhood of a knot $L$ in a closed solid
torus $V$. Suppose that $L$ has winding number $w \neq 0$ in the
solid torus.  Fix an embedding of $V$ in $S^3$ and define
$\partial_0 N$ and  $\partial_1 N$ as in the previous section.
Let $m_0$, $m_1$, $l_0$, and $l_1$ be meridians and longitudes
of $N$, and denote their classes in $H_1(N)$ by $[m_0]$, $[m_1]$,
 $[l_0]$, and $[l_1]$. As shown in Appendix
\ref{homology_calculations_appendix},
$H_1(N)$ has abelian presentation
\begin{eqnarray*}
H_1(N) & = & \Big\langle
[m_0],[m_1],[l_0],[l_1] \Big| [l_1] = w \cdot [m_0],
[l_0] = w \cdot [m_1] \Big\rangle \\
& = & \big\langle [m_0],[m_1] \big\rangle.
\end{eqnarray*}

Let $p$ be a prime not dividing $w$.  Set
$$\bar{\theta}:\pi_1(N) \to H_1(N) / p H_1(N)$$
to be the composition of the Hurewicz homomorphism and the quotient map.
The following diagram commutes:
$$
\xymatrix{
 {H_1(\partial_k N) \cong \pi_1(\partial_k N)} \ar[d]_{i_{k*}}
                  \ar[r]^{\,\,\,\,\,\,\,\,\,\,\,\,\,\,\,i_k^*}
            & {\pi_1(N)} \ar[d]^{\bar{\theta}} \\ 
 {H_1(N)} \ar[r]_{\mathrm{mod} \,p \,\,\,\,\,\,\,\,\,
\,\,\,\,\,\,} & {H_1(N) / p H_1(N)} }$$
$i_{k*}$ is an injection with image $U = \langle [m_k],[l_k] \rangle
= \langle [m_k],w[m_{1-k}] \rangle$, and $V = \mathrm{ker}( \mathrm{mod} p)$
is generated by $p[m_k]$ and $p[m_{1-k}]$.  Note that $V \cap U =
\langle [m_k],pw[m_{1-k}] \rangle$.  This implies that
$\mathrm{ker}((\mathrm{mod} \, p) \circ i_k^*) = \mathrm{ker}
(\bar{\theta} \circ i_k^*)$ is the characteristic
subgroup of index $p^2$ in $\pi_1(\partial_k N)$.  It follows that
the boundary components of the cover $\widetilde{N}$ of $N$
corresponding to $\mathrm{ker} (\bar{\theta}) < \pi_1(N)$ are as prescribed in
the statement of the lemma.  Also,
$\widetilde{N}$ has $\left| H_1(N) / p H_1(N)\right| = p^2$ sheets.
\end{proof}

We now consider which primes cannot divide the nonzero winding number
of a piece in the satellite (JSJ) decomposition of the complement
of $K$.  A bound on the winding number in the lemma above will be
of use.

\begin{Lem}
Suppose $N$ is a piece in the JSJ decomposition of the complement of 
a nontrivial knot in $S^3$ with a $c$-crossing diagram.  Suppose
further that $N$ is the complement of an open regular neighborhood
of a knot in a solid torus with winding number $w$ in the solid torus.
Then $w \leq \frac{c}{2}$.
\label{windinglem}
\end{Lem}

\begin{proof}
Let $K$ be a nontrivial knot in $S^3$ and $N$ be a piece of the JSJ
decomposition of $M=S^3-N(K)$.  Further, suppose that $N$ has two
boundary components.  Let $M'$ be the component of
$M-\partial_0 N$ which is disjoint from  $\partial M$. Then  $M'$ is the 
complement of a knot, $K'$, in $S^3$, and $K$ is a satellite of $K'$.  By 
Theorem 3 of \cite{Schubert} the bridge number, $b'$, of $K'$ must be less 
than or equal to the bridge number, $b$, of $K$.

Consider the knot $K'$.  It has a further satellite knot decomposition
since $N$ has two boundary components.  Let  $K''$ be the companion for
this decomposition, and let $b''$ be its bridge number.  The winding
number of $K'$ in the solid torus is $w$.  Using 
\cite[Theorem 3]{Schubert} we can conclude that
\begin{equation*}
wb'' \leq b'.
\end{equation*}
Clearly $2 \leq b''$ so
\begin{equation*}
2w \leq w b'' \leq b' \leq b.
\end{equation*}
The bridge number of a knot must be less than or equal to the crossing 
number of any projection of the knot so $b \leq c$; hence, we get the 
desired result:
\begin{equation*}
2w \leq c
\end{equation*}
\end{proof}

The covers given by Lemma \ref{homology_lemma} are quite nice in that
they can be made to have few sheets relative to the crossing number
of our original knot.  These are the easy cases because homology
does all the work for us.  Let us consider the types of pieces of the
JSJ decomposition of the knot complement, $M$, which are not covered by
Lemma \ref{homology_lemma}.  Given our restriction on $K$,
there are three remaining 
cases:
\begin{enumerate}
\item $M_i$ is a hyperbolic knot complement.
\item $M_i$ is hyperbolic and the complement of a knot in a 
solid torus with winding number 0 in the solid torus.
\item $M_i$ is a torus knot complement.
\end{enumerate}
We will see in section
\ref{seifert_fibered_pieces_of_a_knot_complement_subsection}
that any Seifert fibered piece of $M$ which is the complement of a
knot in a solid torus satisfies the hypotheses of Lemma
\ref{homology_lemma}.




\section{Hyperbolic Pieces}
\label{hyperbolic_pieces_section}




\subsection{Mahler measure and height}
\label{mahler_measure_and_height_subsection}

In order to address the case in which $M_i$ is a hyperbolic manifold we 
will use the number of tetrahedra in an ideal triangulation of 
$M_i$ to limit certain quantities related to
a representation
of $\pi_1(M_i)$ into $SL_2(\mathbb{C})$.
In the process
we will encounter certain polynomial equations and algebraic numbers.
For numerous reasons the
most natural notions of complexity for polynomials and algebraic numbers are
given by the Mahler measure and height, respectively.  We define these
notions here.

Let $P = P(X_1,X_2, \cdots, X_n)$ be a polynomial with complex 
coefficients. As in \cite{Mignotte} the \emph{Mahler measure}, $M(P)$,
 is given by
\begin{equation*}
M(P) = \exp \left (\int_0^1 \cdots \int_0^1 \log | P(e^{2 \pi i t_1},
\cdots ,e^{2 \pi i t_n})| \, dt_1 \cdots dt_n \right ).
\end{equation*}

As mentioned above, the Mahler measure of a polynomial will be a measure 
of its complexity. Another notion of complexity
which may at first seem more natural is the quadratic norm.
If $P(X_1, \dots, X_n) = \sum a_{j_1 \cdots j_n}X_1^{j_1}
\cdots X_n^{j_n}$ is a polynomial with complex coefficients,
then the {\em quadratic norm} of $P$ is
\begin{equation*}
\| P \| = \sqrt{\sum |a_{j_1 \cdots j_n} |^2}.
\end{equation*}
The following lemma from \cite{Mignotte} relates these two
notions.  Lemma 2.1.7 of \cite{Mignotte} is as follows:
\begin{Lem}
If $P \in \mathbb{C}[X_1, \dots, X_n]$, then $M(P) \leq \| P \|$.
\label{quadraticlem}
\end{Lem}

Lemma 2.1.9 of \cite{Mignotte} relates the Mahler measure and
degree of $P \in \mathbb{C}[X]$ to the size of the coefficients of $P$.
\begin{Lem}
Let $P(X) = c_0 +c_1X + \cdots c_mX^m \in \mathbb{C}[X]$ be a 
polynomial in one variable. Then
$$|c_i| \leq {m \choose i} M(P).$$
In particular $|c_0|,|c_m| \leq M(P).$
\label{coeflem}
\end{Lem}

Let $\alpha \in \mathbb{C}$ be an algebraic number of degree $m$ and 
let $P(X)$ be its minimal polynomial over $\mathbb{Z}$.  Define the 
measure $M(\alpha)$ of $\alpha$ to be
\begin{equation*}
M(\alpha) = M(P).
\end{equation*}
Closely related to the measure of $\alpha$ is its {\em absolute
multiplicative height}, $H(\alpha)$, given by the equation
\begin{equation*}
H(\alpha) = M(\alpha)^{1/m}.
\end{equation*}
At times it is more convenient to consider the {\em absolute
logarithmic height}, $h(\alpha)$, of an algebraic number
$\alpha$ given by
\begin{equation}
h(\alpha) = \log H(\alpha) = \frac{1}{m} \log M(\alpha).
\label{logheight}
\end{equation}
Let $\alpha$ and $\beta$ be algebraic numbers.
We have the following facts found in \cite[Lemma 2A]{Schmidt}
\begin{equation*}
H(\alpha \beta) \leq H(\alpha)H(\beta).
\end{equation*}
\begin{equation*}
H(\alpha + \beta) \leq 2H(\alpha)H(\beta).
\end{equation*}
Equivalently,
\begin{equation}
h(\alpha \beta) \leq h(\alpha) + h(\beta).
\label{height_multiplication_inequality}
\end{equation}
\begin{equation}
h(\alpha + \beta) \leq \log 2 + h (\alpha) + h(\beta).
\label{height_addition_inequality}
\end{equation}

There is a natural notion of height for vectors of
algebraic numbers which is defined in \cite[page 192]{Schmidt}.
For our purposes it will be enough to know that if
$\mbox{\boldmath $\alpha$} = (\alpha_1, \dots, \alpha_k)$ is a vector of
algebraic numbers then for all $i$, $1 \leq i \leq k$
\begin{equation*}
H(\alpha_i) \leq H(\mbox{\boldmath $\alpha$}),
\end{equation*}
and
\begin{equation*}
h(\alpha_i) \leq h(\mbox{\boldmath $\alpha$}).
\end{equation*}

A highly nontrivial result due to Shou-Wu Zhang \cite{Zhang}
is as follows:

\begin{Lem}
Let $P_1, P_2, \cdots, P_n$ be polynomials in the variables $X_1, X_2,
\cdots, X_k$. If $\mbox{\boldmath $\alpha$} = (\alpha_1, \alpha_2,
\cdots, \alpha_k)$ is an isolated solution to the equations $P_i = 0$
then its absolute logarithmic height is bounded as follows:
$$h(\mbox{\boldmath $\alpha$}) \leq A(n)\left(\sum_{i=1}^n \deg
P_i\right)
\left(\left(\sum_{i=1}^n \frac{M(P_i)}{\deg P_i}\right)+B(n)
\log 2\right)$$
where $A(n)=(n^2-n+1)!n^{-2}[(n-1)!]^{-n}$
and $B(n)=(n^3-n^2)(n^2-n+1)^{-1}$.
\label{zhanglem}
\end{Lem}

Finally, a technical result that will be used in the discussion of the
hyperbolic pieces is as follows:

\begin{Lem}
Suppose $A = B_1 B_2 \cdots B_k$ is the product of the $k$ matrices,
$$B_i = \left(
\begin{array}{cc} \beta^i_{11} & \beta^i_{12} \\
\beta^i_{21} & \beta^i_{22} \end{array} \right),$$
with each $\beta^i_{jl}$ an algebraic number.
Suppose further that $h(\beta^i_{jl}) \leq h$ for all
$i,j,l$.  Then for each entry, $a_{ij}$, of $A$ 
$$h(\alpha_{ij})  \leq (2^{k-1}-1)( \log 2 ) + (2^k+2^{k-1}-2) h.$$
\label{matrixlem}
\end{Lem}

\begin{proof}
The proof is a straight-forward induction on $k$.  Clearly the lemma holds 
for the case $k=1$.  Now suppose it is true for $A' = B_1 B_2 \cdots 
B_{k-1}$. Let $A = B_1 B_2 \cdots B_{k-1} B_k = A'B_k$.  Notice that
each entry of $A$ is of the form
$\alpha'_{i1} \beta^k_{1j} + \alpha'_{i2} \beta^k_{2j}$
where $\alpha'_{ij}$ is the $i,j$th entry of $A'$ and
$\beta^k_{ij}$ is the $i,j$th entry of $B_k$. Using
inequalities (\ref{height_multiplication_inequality}) and
(\ref{height_addition_inequality}) we may conclude that
{\setlength\arraycolsep{4pt}
\begin{eqnarray*}
h(\alpha_{ij}) & = & h(\alpha'_{i1} \beta^k_{1j} + \alpha'_{i2}
\beta^k_{2j}) \\
& \leq & \log 2 + h(\alpha'_{i1}\beta^k_{1j}) + h(\alpha'_{i2}
\beta^k_{2j}) \\
& \leq & \log 2 + h(\alpha'_{i1}) +
h(\beta^k_{1j}) + h(\alpha'_{i2}) + h(\beta^k_{2j}) \\
& \leq & \log 2 + 2h + 2 (2^{k-2}-1)( \log 2 ) + 2 (2^{k-1}+2^{k-2}-2)h \\
& = & (2^{k-1}-1)( \log 2 ) + (2^k+2^{k-1}-2) h.
\end{eqnarray*}}
\end{proof}




\subsection{Covering hyperbolic pieces}
\label{covering_hyperbolic_pieces_subsection}

We are now ready to produce the desired covers of the hyperbolic pieces
in the JSJ decomposition of our knot complement which do not
satisfy Lemma \ref{homology_lemma}.  Let $D(n)$ be as in
(\ref{primebound}).

\begin{Thm}
Suppose $N - \partial N$ has a complete, finite volume 
hyperbolic structure, and $N$ is either the complement of a
knot in $S^3$ or the complement of a knot in a solid torus
with winding number $0$.
If $N$ has a combinatorial ideal triangulation with $n$ tetrahedra
then there is a number $B \in \mathbb{N}$ with $B \leq D(n)$
such that for every prime $p$ not dividing $B$,
$N$ has a finite cover $\widetilde{N}$
with at most $p^{3\left( 2^{4n+4} \left(8n^2 + 4n \right)^{
(4n+4)2^{4n+4}}\right)}$ sheets
in which each boundary component of $\widetilde{N}$ is the
noncyclic $p^2$-sheeted cover of a boundary component of $N$.
\label{hypthm}
\end{Thm}

\begin{proof}
Suppose $N - \partial N$ has a complete, finite volume hyperbolic
structure, and is either a piece of the JSJ decomposition of the
complement of a knot in $S^3$ or a knot in a solid torus with winding
number $0$.  Suppose further that $N$ has  a combinatorial ideal
triangulation with $n$ tetrahedra.  Let $p$ be a prime integer.

$N$ has one or two boundary components: $\partial_0 N$ and
possibly $\partial_1 N$.  Ignoring questions of base points
for the moment, we will produce a homomorphism from $\pi_1(N)$
to a finite group whose kernel will intersect each
$\pi_1(\partial_k N)$ in the subgroup
$p \cdot \pi_1(\partial_k N)$.
The covering space of $N$ corresponding to this kernel will
be the desired cover.

More explicitly, let $m_k$ and $l_k$ be a meridian and longitude
for $\partial_k N$.  Choose the base point of $\partial_k N$
to be the point of intersection of $m_k$ and $l_k$, and fix
a path from the base point of $N$ to the base point of $\partial_k N$.
We get explicit inclusions $i_k^*: \pi_1(\partial_k N) \to
\pi_1(N)$.  Let $\lambda_k, \mu_k \in \pi_1(N)$ be the classes
of $m_k$ and $l_k$ respectively. For an oriented loop, $b$,
in the space $U$ let $[b]_U \in H_1(U)$ denote its homology class.
If no space $U$ is indicated then we will assume the space is $N$. 
As in section \ref{covers_from_homology_section} set the
homomorphism
$$\bar{\theta}:\pi_1(N) \to H_1(N) / pH_1(N)$$ to be the composition
of the Hurewicz map and the quotient map.

The manifold, $N$, is either the complement of a knot in $S^3$ or
the complement of a knot with winding number $0$ in the solid torus.
In either case we have that $\bar{\theta}(\lambda_k) = 0$ and
$\bar{\theta}(\mu_k)$ has order exactly $p$.
Using the hyperbolic structure of $N$ we will produce another
homomorphism
$$\bar{\rho}:\pi_1(N) \to SL(2,F)$$ for $F$ some finite field of
characteristic $p$.  By construction it will be clear that
 $\bar{\rho}(\lambda_k)$
has order exactly $p$, and $\bar{\rho}(\mu_k)$ has order dividing
$p$.  It follows directly that
$\bar{\rho} \times \bar{\theta}$,  will
be a homomorphism with the desired kernel.  The challenge will
be to show that for all $p > D(n)$ such a $\bar{\rho}$ exists
and to bound the minimum degree of the finite field $F$ over
$\mathbb{Z} / p \mathbb{Z}$ from above. Clearly, this will
bound the order of the group $SL(2,F) \times \mathbb{Z} / p
\mathbb{Z} \times \mathbb{Z} / p \mathbb{Z}$ which will
contain the image of $\bar{\rho} \times \bar{\theta}$.

Produce a presentation of $\pi_1(N)$ as follows:
Let $C$ be the standard spine of $N$ dual to $\mathcal{T}$.
Clearly, $C$ will have $n$ vertices (one for each ideal tetrahedron of
$\mathcal{T}$) and $4n/2 = 2n$ edges (one for each face of
$\mathcal{T}$).  Note that $C$ is homotopy equivalent to $N$,
and $\partial N$ is a union of tori.  It follows that $C$ has
Euler characteristic $0$.  This implies that $C$ must have $n$
faces.  If we fix a maximal tree in the $1$-skeleton of $C$, we
get a presentation $\langle g_0, g_1, \dots, g_n | r_1, r_2, \dots
, r_n \rangle$ for $\pi_1(C) \cong \pi_1 (N)$ with $n+1$ generators
and $n$ relations.  Furthermore, each edge of $C$ is incident with 3
faces so the sum of the lengths of the relations must be $3(n+1)$.

The ideal triangulation, $\mathcal{T}$, of $N$ with $n$ ideal tetrahedra
induces a natural triangulation of $\partial N$ as follows:
Place a single normal triangle in each corner of each ideal
tetrahedron of $\mathcal{T}$. Gluing these triangles to form a normal
surface gives a triangulation for $\partial N$ with exactly $4n$
triangles.  Dual to this triangulation is a polygonal decomposition
of the boundary tori with $4n$ vertices whose $1$-skeleton is a
trivalent graph.  For each boundary component, $\partial_k N$, we have
paths $x_k'$, and $y_k'$ in this $1$-skeleton generating the
fundamental group of that boundary component.  In fact we may
assume that these paths all have length at most $4n$ by insisting
that $x_k'$ and $y_k'$ traverse each vertex at most once.  These
paths project in a natural way onto the $1$-skeleton of the
standard spine dual to $\mathcal{T}$. Whence we get words
$x_k$ and $y_k$ in $g_0, g_1, \dots, g_n$
generating the fundamental group of $\partial_k N$ as a subgroup of
$\pi_1 (N)$.  Furthermore, the words $x_k$ and $y_k$ have length
at most $4n$.

In order control the image of $\lambda_k$ under $\bar{\rho}$, we
will bound its length as a word in $x_k$ and $y_k.$ (We
may assume after adjusting paths connecting
base points that $\lambda_k \in \langle
 x_k, y_k \rangle$.)

Assume for the moment that $N$ is a knot complement. Then
the homology class $[l_0]$ is trivial.
Consider the presentation
$\langle g_0, g_1, \dots, g_n | r_1, r_2, \dots
, r_n \rangle$ for $\pi_1(N)$.  There is a corresponding abelian
presentation for $H_1(N)$.
For $0 \leq j \leq n$ choose
$\nu_j \in \mathbb{Z}$ so that $[g_j] = \nu_j [m_0].$
Each relation,
$r_j$, is trivial in $\pi_1(N)$ and hence must map to $0$ in
$H_1(N)$.  This
translates to an equation specifying that some integral linear
combination of $\nu_j$'s is $0$.  For example, the relation
$g_0 g_2 g_3 g_1^{-1} g_2 = 1$ would give the equation $1\nu_0 -1\nu_1 + 
2\nu_2 + 1\nu_3 = 0$. Consider the $n \times (n+1)$ matrix, $B$, whose
$ij$th entry is the coefficient of $\nu_j$ in the
equation coming from the relation $r_i$. The vector
$(\nu_0, \cdots, \nu_n)$ will be the smallest
nonzero, integral vector whose dot product with each row of $B$ is
$0$. Since $B$ is a presentation matrix for the homology
of $N$, this property actually characterizes $(\nu_0, \cdots, \nu_n)$
up to sign. Let $B_j$ be the $n \times n$ minor of $B$
formed by dropping
the $j$th column of $B$.  I claim that the integer vector
$(\det B_0, -\det B_1, \det B_2, \cdots, (-1)^n \det B_n)$
is a multiple of $(\nu_0, \cdots, \nu_n)$.  Let $\mathbf{w} =
(w_0, \cdots, w_n)$ be an arbitrary vector.  By definition,
$\mathbf{w}$ will be in the
row space of $B$ if and only if the determinant of the $n \times n$
matrix formed by adding $\mathbf{w}$
in as the first row of $B$ is $0$.  Hence,
$\mathbf{w}$ will be in the row space of $B$ if and only if
$\sum_{j=0}^n w_j (-1)^j \det (B_j) = 0$.  This shows that
$(\det B_0, -\det B_1, \det B_2, \cdots, (-1)^n \det B_n)$ is
indeed perpendicular to the row space of $B$ and must be a multiple
of $(\nu_0, \cdots, \nu_n)$.

Recall that in a standard spine, an edge in the $1$-skeleton is incident
with exactly $3$ faces.  This implies that the sum of the absolute values
of entries in a column of $B$ is $3$.  It follows that
$|\det (B_j)| \leq 3^n$ for each minor $B_j$.  We have also shown that
$|\nu_j| \leq |\det (B_j)|$ for all $j$, thus $|\nu_j| \leq 3^n$ for all
$j$.

The words $x_0$ and $y_0$ are of length at most
$4n$ in the $g_i$'s.  Hence if $a_0,b_0 \in \mathbb{Z}$ are
chosen so that $[x_0] = a_0 [m_0]$ and $[y_0] = b_0 [m_0]$,
then $|a_0|,|b_0| \leq 4n3^n$.  In
$H_1(\partial_0 N)$ we can write the homology class of $l_0$
as a linear combination of $[x_0]_{\partial_0 N}$ and
$[y_0]_{\partial_0 N}$ by noting that if $[l_0]_{\partial_0 N}
= v [x_0]_{\partial_0 N} + w [y_0]_{\partial_0 N}$ then $(v,w)$ generates
the null space of the $1 \times 2$ matrix $(a_0 \,\, b_0)$.  In fact
\begin{equation}
\pm [l_0]_{\partial_0 N} = -b_0 [x_0]_{\partial_0 N} + a_0
[y_0]_{\partial_0 N}
\label{longitude_length_equation}
\end{equation}
where $|a_0|,|b_0| \leq 4n3^n$.

Now suppose $N$ is the complement of a knot with winding number $0$
in the solid torus.  Then $H_1(N)$ has abelian presentation
$$H_1(N) = \Big\langle
[m_0],[m_1],[l_0],[l_1] \Big| [l_1] = 0,
[l_0] = 0 \Big\rangle$$
As above, let $B$ be the presentation matrix for $H_1(N)$ coming
from the presentation of
$\pi_1(N)$.  For each $k \in \{0,1\}$ at least one of $[x_k]$ or
$[y_k]$ must be nontrivial.  Without
loss of generality assume $[x_k]$ is nontrivial.  Let $B^{(k)}$
be the $(n+1) \times(n+1)$ matrix whose first $n$ rows agree with
 $B$ and whose $(n+1)$th row comes from the word $x_k$.  Let
$\nu_0^{(k)}, \cdots, \nu_n^{(k)} \in \mathbb{Z}$ be integers
such that
$$[g_i] = \nu_i^{(0)}\cdot [m_0] + \nu_i^{(1)}\cdot [m_1]$$
Then $(\nu_0^{(k)}, \cdots, \nu_n^{(k)})$ generates the null
space of $B^{(1-k)}$.  One of the first $n$ rows of $B$ is
a linear combination of the others (over $\mathbb{Q}$).  If
we remove this row, $(\nu_0^{(k)}, \cdots, \nu_n^{(k)})$
still generates the null space of the new matrix.  The same
argument as above shows that $|\nu_j^{(k)}| \leq 16n^2 3^{n-1}$.
We may then conclude that there are integers $a_k,b_k \in \mathbb{Z}$
such that 
\begin{equation}
\pm [l_k]_{\partial_k N} = -b_k [x_k]_{\partial_k N} + a_k [y_k]_{\partial_k N}
\label{two_longitude_length_equation}
\end{equation}
where
$$|a_k|,|b_k| \leq 16n^2 3^{n-1}.$$

The hyperbolic structure of $N$ gives a faithful representation
of $\pi_1(N)$ in $PSL_2(\mathbb{C})$.  As Thurston has shown, we
may lift this representation
to a faithful representation $\rho: \pi_1(N) \to SL_2(\mathbb{C})$.
Let
$$\rho(g_k) = \left( \begin{array}{cc} z_{4k} & z_{4k+1} \\
z_{4k+2} & z_{4k+3} \end{array} \right)$$
Let $A = \mathbb{Z}[z_0,z_1,\cdots,z_{4n+3}] \in \mathbb{C}$.
Since
$$\rho(g_k^{-1}) = \left( \begin{array}{cc} z_{4k+3} & -z_{4k+2} \\
-z_{4k+1} & z_{4k} \end{array} \right),$$
we actually have $\rho: \pi_1(N) \to SL_2(A)$.  If we produce a
ring homomorphism $\eta:A \to F$ for some finite field $F$ then
it will induce a group homomorphism $\eta*:SL_2(A) \to SL_2(F)$.
Composing $\eta *$ with $\rho$ will give $\bar{\rho}:\pi_1(N) \to SL(2,F)$.
Our goal is to bound the height and degree of $\mathbf{z} =
(z_0,z_1,\cdots,z_{4n+3}) \in \mathbb{C}^{4n+4}$, and
extract a sufficient criterion on $F$ for it to have a
suitable nontrivial ring homomorphism $\eta:A \to F$.

The point $\mathbf{z} \in \mathbb{C}^{4n+4}$ satisfies
the $n+1$ polynomial equations
$$G_k = X_{4k} X_{4k+3} - X_{4k+2} X_{4k+1} - 1 = 0$$
specifying that $\det \rho (g_k) = 1$.

Each relation $r_k$
gives four polynomial relations satisfied
by $\mathbf{z}$ indicating that $\rho (r_k)$ is
the identity matrix.  One of these equations is superfluous, so we
drop it and let $R_{3k}$, $R_{3k+1}$, and $R_{3k+2}$ be the remaining three.
If the relation $\rho_k$ has length
$l$ then the polynomials $R_{3k}$,$R_{3k+1}$, and $R_{3k+2}$ will have
degree at most $l$ and will be the sum of at most $2^l$ monomials
with coefficient $\pm 1$ and possibly the term $-1$.

After an appropriate conjugation we may assume without loss
of generality that
$$\rho(x_0) = \left( \begin{array}{cc} 1 & 1 \\ 0 & 1 \end{array}
\right).$$
Hence we get four more polynomial relations satisfied by $\mathbf{z}$.
Again, we drop one of them and let $C_0$, $C_1$, and $C_2$ be the three
that remain.  By the same argument
as above each polynomial relation $C_i$ is a sum of at most $2^{4n}$
terms of degree at most $4n$ and possibly the term $-1$.

If $N$ has two boundary components then by the completeness of
$N$ we must have $\mathrm{tr} (\rho(x_1)) = \pm 2$. We
specify this with the single polynomial relation $Q$.
The same argument gives that $Q$ is a sum of at most $2^{4n}$
terms of degree at most $4n$ and the term $\mp 2$.

In fact, $\mathbf{z} \in \mathbb{C}^{4n+4}$ is an isolated
root of the $4n+4$ or $4n+5$ polynomials, $\mathcal{P} = \{ G_0,\dots,G_n$,
$R_3,\dots,R_{3n+2}$, $C_0,C_1,C_2,Q \}$.
(See \cite[Proposition 2]{Culler}.)
In light of Lemma \ref{zhanglem}, bounds on the Mahler measures
of these polynomials will give a bound on the height of $\mathbf{z}$.

We now bound the Mahler measures of these polynomials.
Firstly by Lemma \ref{quadraticlem},
$$ M(G_k) = M(XY-ZW-1) \leq \| XY-ZW-1 \| = \sqrt{3}$$
If the relation $r_k$ has length $l$ then polynomials $R_{3k},
 R_{3k+1},$ and $R_{3k+2}$ will be sums of  $2^l$ monomials
with coefficient $\pm 1$ and possibly the term $-1$.
Thus if $\mathbf{w}$ is a vector of complex numbers with norm 1
then clearly $| R_i(\mathbf{w}) | < 2^l + 1$.  We get the following bound.
\begin{eqnarray*}
M(R_i) & = & \exp \left (\int_0^1 \cdots \int_0^1 \log | R_i(e^{2 \pi i t_1},
\cdots ,e^{2 \pi i t_{4n+4}})| \, dt_1 \cdots dt_{4n+4} \right ) \\
& \leq & \exp \left (\int_0^1 \cdots \int_0^1 \log (2^l+1) \, dt_1 \cdots dt_{4n+4} \right ) \\
& = & 2^l+1.
\end{eqnarray*}
Similarly $M(C_i) \leq 2^{4n} + 1$, and $M(Q) \leq 2^{4n} + 2$.

We now have all the necessary ingredients to bound $h( \mathbf{z} )$.
The degree of each $G_i$ is 2.  The sum of the lengths of the
relations $r_i$ is $3(n+1)$, so the sum of the degrees of the
$R_i$'s is at most $9(n+1)$.  Each $C_i$ has degree at most $4n$,
and $Q$ has degree at most $4n$ as well. By Lemma \ref{zhanglem},

\begin{eqnarray}
h(z_k) & \leq & h(\mathbf{z}) \nonumber \\
& \leq &  A(4n+5)\left( \sum_{P \in \mathcal{P}} \deg P \right)
\left(\left(\sum_{P \in \mathcal{P}} \frac{M(P)}{\deg P}\right)
+B(4n+5)\log 2\right) \nonumber \\
& \leq & A(4n+5)\left( 2(n+1) + 3 \cdot 3(n+1) + 3\cdot 4n + 4n \right)
\label{z_height_bound_inequality} \\
& & \cdot \left( \left( {\scriptstyle \frac{\sqrt{3}}{2}}(n+1)
+ 3 \cdot \left(2^{3(n+1)}+n \right) + 4 \cdot 2^{4n} \right)
+B(4n+5)\log 2\right) \nonumber \\
& \leq & A(4n+5)( 27n + 5) \nonumber \\
& & \cdot \left( 2^{4n+2} + 3\cdot 2^{3n+3} + \left({\scriptstyle 
\frac{\sqrt{3}}{2}} + 3 \right) n + B(4n+5)\log 2
+{\scriptstyle \frac{\sqrt{3}}{2}} \right). \nonumber
\end{eqnarray}

The ring homomorphism, $\eta :A \to F$, is given as follows.  Let
$W_k(X) \in \mathbb{Z} [X]$ be the minimal polynomial
of $z_k$.
Fix a prime, $p$, and let $W_k^{(p)}(X)$ denote the image of $W_k(X)$
under the natural map $\mathbb{Z} [X] \to \mathbb{F}_p [X]$
where $\mathbb{F}_p = \mathbb{Z} / p\mathbb{Z}.$
If $p$ does not divide the leading coefficient of $W_1(X)$
then we are assured that $W_1^{(p)}(X)$ has a root, $\zeta_1$,
in the algebraic closure of $\mathbb{F}_p$, and we have a ring
homomorphism $\mathbb{Z}[z_1] \to \mathbb{F}_p(\zeta_1)$ taking
$z_1$ to $\zeta_1$. Set $A_1 = \mathbb{Z}[z_1]$ 
and $F_1 = \mathbb{F}_p(\zeta_1)$ In this case we get a
ring homomorphism, $\eta_1: A_1 \to F_1$.  Let $\tilde{\eta}_1:A_1[X]
\to F_1[X]$ be the induced map on the polynomial rings.

Inductively let $S_{i+1}(X) \in A_i[X]$ be the minimal polynomial of $z_{i+1}$
over $A_i$.  If $p$ does not divide the leading coefficient of $W_{i+1}(X)$
then $W_{i+1}^{(p)}(X)$ has roots in the algebraic
closure of $\mathbb{F}_p$ some of which will be roots of
$\tilde{\eta}_i\left( S_{i+1}(X) \right)$. Let $\zeta_{i+1}$ be one
such root.
Set $A_{i+1} = A_i[z_{i+1}]$ 
and $F_{i+1} = F_i(\zeta_{i+1})$. In this case we get a
ring homomorphism $\eta_{i+1}: A_{i+1} \to F_{i+1}$ restricting to
$\eta_i$ on $A_i$ and taking $z_{i+1}$ to $\zeta_{i+1}$.
Let $\tilde{\eta}_{i+1}:A_{i+1}[X] \to F_{i+1}[X]$ be the induced map on the
polynomial rings.

From this discussion it is clear that we will have a homomorphism,
$\eta_{4n+4}:A_{4n+4} \to F_{4n+4}$, if $p$ does not divide any of the leading
coefficients of the $W_k$'s.  Of course, $A = A_{4n+4}$.
Set $\eta = \eta_{4n+4}$ and $F = F_{4n+4}$.
After further restriction on $p$, $\eta$ will be the desired
homomorphism.

We will now proceed to bound the degree $[F:\mathbb{F}_p]$.
In \cite{Dube} a bound on the degree of a polynomial in
a Gr\"{o}bner basis with any monomial order is given.  This gives a bound
on the degrees of the polynomials in a Gr\"{o}bner basis of the ideal generated
by $\mathcal{P}$ which in turn gives a bound on the
degree of $W_k$.
\begin{eqnarray*}
\deg (W_k) & \leq & 2(8n^2+4n)^{2^{4n+4}}
\label{w_degree_bound_inequality}
\end{eqnarray*}
Thus, 
\begin{eqnarray*}
[F:\mathbb{F}_p] & =  & \prod_{i=1} \deg \big[\tilde{\eta}_i\left( S_{i+1}(X)
\right) \big] \\
& \leq & \prod_{i=1} \deg W_i \\
& \leq & 2^{4n+4} \left( 8n^2 + 4n \right)^{(4n+4)2^{4n+4}}
\end{eqnarray*}
We also get a bound of the order of the field $F$.
\begin{equation*}
|F| \leq p^{\left( 2^{4n+4} \left( 8n^2 + 4n \right)^{(4n+4)2^{4n+4}}\right) }.
\end{equation*}
The order of $SL_2(F)$ is $(|F|^2 - 1)(|F| - 1)$; hence,
\begin{equation}
|SL_2(F)| \leq |F|^3 \leq p^{3\left( 2^{4n+4} \left(
8n^2 + 4n \right)^{(4n+4)2^{4n+4}}\right) }.
\label{slsize}
\end{equation}

We now address the question of how to ensure that $p$ meets
all of the conditions stipulated above.
Inequalities (\ref{z_height_bound_inequality}) and
(\ref{w_degree_bound_inequality}) combine to bound the Mahler
measure of $z_k$.
\begin{eqnarray}
M(z_k) & \leq &\exp \Big[ 2\left( 8n^2 + 4n \right)^{2^{4n+4}}
A(4n+5)( 27n + 5) \nonumber \\ [-3mm]
& & \label{z_measure_bound_inequality} \\ [-3mm]
& & \cdot \left( 2^{4n+2} + 3\cdot 2^{3n+3} + \left({\scriptstyle 
\frac{\sqrt{3}}{2}} + 3 \right) n +B(4n+5)\log 2
+{\scriptstyle \frac{\sqrt{3}}{2}} \right) \Big]. \nonumber
\end{eqnarray}
If $p$ does not divide the coefficient the highest degree term
of $W_k$ then $S_k$ has roots in $\Omega_p$.  Lemma \ref{coeflem}
shows that the coefficient of the highest degree term of $W_k$ is
less than $M(W_k)$. This is sufficient to ensure that
we have a ring homomorphism $\eta :A \to F$.  However, we need further
restrictions on $p$ to ensure that $\rho (l_i)$ is nontrivial.

Recall that $[l_0]_{\partial_0 N} = -b_0[x_0]_{\partial_0 N} + a_0[y_0]_{\partial_0 N}$ for some $a_0,b_0 \in \mathbb{Z}$ with
$|a_0|,|b_0| \leq 16n^2 3^{n-1}$.
$$\rho(x_0) = \left( \begin{array}{cc} 1 & 1 \\
0 & 1 \end{array} \right)$$
and
$$\rho(y_0) = \left( \begin{array}{cc} 1 & \alpha_0 \\
0 & 1 \end{array} \right)$$
The word $y_0$ has length at most $4n$ so by Lemma \ref{matrixlem}
\begin{eqnarray*}
h(\alpha_0) & \leq & (2^{4n-1}-1)( \log 2 ) + (2^{4n}+2^{4n-1}-2)
A(4n+5)( 27n + 5) \\
& & \cdot \left( 2^{4n+2} + 3\cdot 2^{3n+3} + \left({\scriptstyle 
\frac{\sqrt{3}}{2}} + 3 \right) n 
+B(4n+5)\log 2  +{\scriptstyle \frac{\sqrt{3}}{2}} \right).
\end{eqnarray*}
It follows that
$$\rho(l_0) = \left( \begin{array}{cc} 1 & -b_0+a_0 \alpha_0 \\
0 & 1 \end{array} \right)$$
and
\begin{eqnarray*}
h(-b_0+a_0 \alpha_0) & \leq & \log 2 + \log a_0 + \log b_0 + h(\alpha) \\
& \leq & \log 2 + 2 \log (16n^2 3^{n-1}) + (2^{4n-1}-1)( \log 2 ) \\
& & + (2^{4n}+2^{4n-1}-2)A(4n+5)( 27n + 5) \\
& & \cdot \left( 2^{4n+2} + 3\cdot 2^{3n+3} + \left({\scriptstyle 
\frac{\sqrt{3}}{2}} + 3 \right) n 
+B(4n+5)\log 2  +{\scriptstyle \frac{\sqrt{3}}{2}} \right)
\end{eqnarray*}
The degree of $-b_0+a_0 \alpha_0$ is at most the product of the
degrees of the $z_k$'s; hence,
\begin{equation*}
\deg ( -b_0+a_0 \alpha_0 ) \leq 2^{4n+4} \left( 8n^2 + 4n \right)^{(4n+4)2^{4n+4}}.
\end{equation*} 
From (\ref{logheight}) we get:
\begin{eqnarray}
\lefteqn{M(-b_0+a_0 \alpha_0 ) \leq} \nonumber \\
& & \exp
        \Big[ 2^{4n+4} \left( 8n^2 + 4n \right)^{(4n+4)2^{4n+4}}
                \Big(
\log 2 + 2 \log (16n^2 3^{n-1}) \nonumber \\
& & + (2^{4n-1}-1)( \log 2 ) + (2^{4n}+2^{4n-1}-2)A(4n+5)( 27n + 5) 
\nonumber \\ [-3mm]
& & \label{longitude_measure_bound_inequality} \\ [-3mm]
& & \cdot \left( 2^{4n+2} + 3\cdot 2^{3n+3} + \left({\scriptstyle 
\frac{\sqrt{3}}{2}} + 3 \right) n 
+B(4n+5)\log 2  +{\scriptstyle \frac{\sqrt{3}}{2}} \right) \Big) \Big].
        \nonumber
\end{eqnarray}

If $p$ does not divide the constant term of the minimal polynomial
of $-b_0+a_0 \alpha_0$ over $\mathbb{Z}$ then $\eta ( -b_0+a_0 \alpha_0) $
cannot be 0.  The order of $\bar{\rho} (l_0)$ is the additive
order of  $\eta ( -b_0+a_0 \alpha_0) $ which must be $p$.

In the case that $N$ has two boundary components, we also require
 $\bar{\rho} (l_1)$ to have order exactly $p$.  Here we note that
$\rho (x_1)$ is parabolic and so
$$
\rho (x_1) = \left( \begin{array}{cc} a & -b \\
-c & d \end{array} \right)
\left( \begin{array}{cc} 1 & 1 \\
0 & 1 \end{array} \right)
\left( \begin{array}{cc} a & -b \\
-c & d \end{array} \right)^{-1}
=\left( \begin{array}{cc} 1+ac & a^2 \\
-c^2 & 1-ac \end{array} \right)
$$
for some $a,c \in \mathbb{C}$.  Adjoin square roots of the
upper right and lower left entries of $\rho (x_1)$ to $A$ get
the ring $A'$, and extend $\eta:A \to F$ to some ring
homomorphism $\eta':A' \to F'$ (no further restriction on $p$ is needed
for such an $\eta'$ to exist.)
Then in $SL_2(F')$, $\bar{\rho} (x_1)$ will be conjugate to
{\tiny $\left( \begin{array}{cc} 1 & 1 \\
0 & 1 \end{array} \right)$}
as long as at least one of $\eta (a^2)$ or
$\eta (c^2)$ is nonzero.  At least one of $a^2$ and $c^2$ is
nonzero in $\mathbb{C}.$  Assume, without loss of
generality, that $a^2 \neq 0$.  The word
$x_1$ has length at most $4n$.  As above this
gives a bound on the height and degree of $a^2$ which gives
the following bound on the Mahler measure of $a^2$
\begin{eqnarray}
\lefteqn{M(a^2) \leq} \nonumber \\
& &  \exp
        \Big[ 2^{4n+4} \left( 8n^2 + 4n \right)^{(4n+4)2^{4n+4}}
                \Big( (2^{4n-1}-1)( \log 2 ) \nonumber \\
& & + (2^{4n}+2^{4n-1}-2)A(4n+5)( 27n + 5) \nonumber \\ [-3mm]
& & \label{cojugacy_measure_bound_inequality} \\ [-3mm]
& & \cdot \left( 2^{4n+2} + 3\cdot 2^{3n+3} + \left({\scriptstyle 
\frac{\sqrt{3}}{2}} + 3 \right) n 
+B(4n+5)\log 2  +{\scriptstyle \frac{\sqrt{3}}{2}} \right) \Big) \Big].
        \nonumber
\end{eqnarray}
If $p$ does not divide the constant term of the minimal polynomial
of $a^2$ over $\mathbb{Z}$ then $\eta (a^2)$ cannot be 0.  The
above bound on the Mahler measure of $a^2$ is also a bound on the
constant term in the minimal polynomial of $a^2$.

Now consider $\rho (y_1)$. If $a,b,c,d$ are as above, we must have
$\alpha_1 \in \mathbb{C}$ such that
$$
\rho (y_1) = \left( \begin{array}{cc} a & -b \\
-c & d \end{array} \right)
\left( \begin{array}{cc} 1 & \alpha_1 \\
0 & 1 \end{array} \right)
\left( \begin{array}{cc} a & -b \\
-c & d \end{array} \right)^{-1}
$$
If we reverse the labels on the boundary components for a moment
it is clear that the height and degree of $\alpha_1$ satisfy
the same bounds as the ones given for height and degree
of $\alpha_0$.  Again we have $a_1,b_1 \in \mathbb{Z}$ with
$|a_1|,|b_1| \leq 16 n^2 s^{n-1}$ such that $[l_1]_{\partial_1 N}
= -b_1 [x_1]_{\partial_1 N} + a_1 [y_1]_{\partial_1 N}$.
If follows that
$$
\rho (l_1) = \left( \begin{array}{cc} a & -b \\
-c & d \end{array} \right)
\left( \begin{array}{cc} 1 & -b_1 + a_1 \alpha_1 \\
0 & 1 \end{array} \right)
\left( \begin{array}{cc} a & -b \\
-c & d \end{array} \right)^{-1}
$$
$M( -b_1 + a_1\alpha_1)$ satisfies the same bound as was given for
$M( -b_0 + a_0\alpha_0)$.  If $p$ does not divide the constant term of
the minimal polynomial of $-b_1 + a_1 \alpha_1$, then $\bar{\rho} (l_1)$
has order $p$.

In summary, if $p$ does not divide the top degree terms of
the minimal polynomials of the $z_i$'s with $1 \leq i \leq 4n+4$
then $\eta:A \to F$ exists.  If $p$ does not divide
the constant terms of minimal polynomials of $-b_0 + a_0 \alpha_0$,
$a^2$, and $-b_1 + a_1 \alpha_1$ then $\bar{\rho} (l_i)$ has
order exactly $p$.  Let $B$ be the product of all these coefficients.
From Lemma \ref{coeflem} and inequalities
(\ref{z_measure_bound_inequality}),
(\ref{longitude_measure_bound_inequality}), and
(\ref{cojugacy_measure_bound_inequality})
it follows that
\begin{eqnarray}
B & \leq & \exp
        \Big[ 2(4n+4)\left( 8n^2 + 4n \right)^{2^{4n+4}}
A(4n+5)( 27n + 5) \nonumber \\
& & \cdot \left( 2^{4n+2} + 3\cdot 2^{3n+3} + \left({\scriptstyle 
\frac{\sqrt{3}}{2}} + 3 \right) n 
+B(4n+5)\log 2  +{\scriptstyle \frac{\sqrt{3}}{2}} \right) \nonumber \\
& & + 2^{4n+4} \left( 8n^2 + 4n \right)^{(4n+4)2^{4n+4}}
         \Big( 2\log 2 + 4 \log (16n^2 3^{n-1}) \nonumber \\
& & + 3(2^{4n-1}-1)( \log 2 ) + 3(2^{4n}+2^{4n-1}-2)A(4n+5)( 27n + 5) 
\nonumber \\
& & \cdot \left( 2^{4n+2} + 3\cdot 2^{3n+3} + \left({\scriptstyle 
\frac{\sqrt{3}}{2}} + 3 \right) n 
+B(4n+5)\log 2  +{\scriptstyle \frac{\sqrt{3}}{2}} \right) \Big) \Big].
        \nonumber \\
& = & D(n). \nonumber
\end{eqnarray}
If $p$ does not divide $B$ then $\bar{\rho} \times \bar{\theta}$
has the desired kernel.
\end{proof}




\section{Seifert Fibered Pieces}
\label{seifert_fibered_pieces_section}




\subsection{Seifert fibered pieces of a knot complement}
\label{seifert_fibered_pieces_of_a_knot_complement_subsection}

We now turn our attention to Seifert fibered pieces of our knot 
complement.
\begin{Lem}
Suppose $N$ is a piece in the JSJ decomposition of the complement of
a nontrivial, decompositionally linear knot 
in $S^3$, and $N$ is Seifert fibered.  Then $N$ has base orbifold either a 
disk with two cone points of order $u$ and $v$ with $(u,v)=1$, or $N$ 
has base orbifold an annulus with one cone point.
\label{sieftypelem}
\end{Lem}

\begin{proof}
As in the statement of the lemma we assume $N$ is a Seifert
fibered piece in the JSJ
decomposition of the complement of a decompositionally linear knot in $S^3$.
  Then $N$ has one or two boundary components.

\emph{Case 1.} Assume $N$ has one boundary component.  Then since
$N$ is embedded
in $S^3$, the boundary of $N$ is a torus in $S^3$.  By the Solid
Torus Theorem (see \cite[page 107]{Rolfsen}), 
$\partial N \subset S^3$ bounds a solid torus on at least one side.
Clearly, $N$ cannot be a solid torus (since
we assumed the knot to be nontrivial), so
$S^3 - N$ must be a solid torus. This shows that $N$ is actually
the complement of a knot, $K'$, in $S^3$.  The only knots with
Seifert fibered complements are torus knots (See
\cite[Theorem 10.5.1]{Kawauchi}). Hence,
for some relatively prime $u,v \in \mathbb{Z}$, the knot $K'$
is the $(u,v)$-torus knot.  It follows that $N$ has a Seifert
fibered structure with two singular fibers with orders $u$ 
and $v$.  Thus, the base orbifold of $N$ is a disk with two cone
points with relatively prime orders $u$ and $v$.
 
\emph{Case 2.} Assume $N$ has two boundary components. The
manifold, $N$ is {\em simple} in the sense that it contains
no essential tori (otherwise we would have cut along them).
As in \cite[Proposition C.5.2]{Kawauchi}, the only
simple Seifert fibered manifolds with two boundary components
have base orbifold an annulus with a single cone point.
\end{proof}

Note that in the case where $N$ has two boundary components,
it is the complement of a knot in the solid torus with winding
number equal to the order of the cone point of its base space.
This situation is covered by Lemma \ref{homology_lemma}.




\subsection{Seifert fibered pieces with one boundary component}
\label{seifert_fibered_pieces_with_one_boundary_component_subsection}

Now one case remains.  If our knot $K$ is a satellite of the 
$(u,v)$-torus knot then the JSJ decomposition of its complement has a 
Seifert fibered piece whose base space is a disk with two cone points. 
We will bound $u$ and $v$ based on the crossing number of $K$.  The 
bridge number of $K$ must be greater than the bridge number of the 
$(u,v)$-torus knot.  The bridge number of a torus knot is known to be 
the smaller of $|u|$ and $|v|$ (See \cite[Theorem 7.5.3]{Murasugi});
consequently, if $c$ is the number of 
crossings in some diagram of $K$ then the smaller of $|u|$ and $|v|$ must 
be less than $c$.  Unfortunately, we must bound the larger of the two.

Recall that Lemma \ref{decomposition_triangulation_lemma} gives a bound on the number of
tetrahedra needed to triangulate any piece in the JSJ decomposition
of the complement of the knot $K$. We will proceed to bound the
minimum number of tetrahedra needed to triangulate the complement
of a $(u,v)$-torus knot from below.  This will be done by showing that a knot
complement that can be triangulated with $n$ ideal tetrahedra has an
Alexander polynomial with degree at most $(n^2+n)3^{n+1}$.

\begin{Lem}
Suppose $L$ is a knot in $S^3$, and its complement $N = S^3 - L$ can
be triangulated with $n$ ideal tetrahedra.  Then the Alexander polynomial
of $L$ has degree at most $(n^2+n)3^{n+1}$.
\label{alexlem}
\end{Lem}

\begin{proof}
Let $L$ be a knot in $S^3$ and $N = S^3 - L$ its complement.  Suppose 
$N$ has an ideal triangulation $\mathcal{T}$ with $n$ ideal tetrahedra.
Set $G = \pi_1(N)$.
As in the proof of Theorem \ref{hypthm}, we have a presentation
$\langle g_0, g_1, \dots, g_n | r_1, r_2, \dots , r_n
\rangle$
for $G$ with $n+1$ generators and $n$
relations.  Furthermore, each edge of $C$ is incident with 3
faces so the sum of the lengths of the relations must be $3(n+1)$.

Following the technique given in \cite[example 9.15]{Burde}
this presentation of the group may be used to find the
first elementary ideal of the Alexander module of $K$.
To begin we let $F = \langle g_0, g_1,
\dots, g_n \rangle$ be the free group.  Let the derivations
$\frac{\partial}{\partial g_i}:\mathbb{Z}F \to \mathbb{Z}F$
be the linear maps satisfying the following rules
for all $\alpha, \beta \in F$:
\begin{itemize}
\item $\frac{\partial}{\partial g_i}(g_j) = \delta_{ij}$.
\item $\frac{\partial}{\partial g_i}(\alpha^{-1}) = 
\alpha^{-1}\frac{\partial}{\partial g_i}(\alpha)$.
\item $\frac{\partial}{\partial g_i}(\alpha \cdot \beta) = 
\frac{\partial}{\partial g_i}(\alpha) + \alpha
\frac{\partial}{\partial g_i}( \beta)$.
\end{itemize}
For each generator $g_i$ and each relation $r_j$ compute
$\frac{\partial}{\partial g_i}r_j$.  Let $\psi:\mathbb{Z}F \to
\mathbb{Z}G$ be the linear extension of the quotient
homomorphism, $F \to G$, and $\varphi: \mathbb{Z}G \to \mathbb{Z}
\langle t \rangle$
be the linear extension of the Hurewicz homomorphism,
$G \to H_1(N) = \langle t \rangle$.
By \cite[Proposition 9.14]{Burde} we know that the ideal
of $\mathbb{Z} \langle t \rangle$ generated by the determinants of the
$n \times n$ minors of $A = \left ( \varphi \circ
\psi (\frac{\partial r_j}{\partial g_i} ) \right )$
will be a principal ideal generated by the Alexander polynomial,
$\Delta(t)$, of $L$.  Let $A_i$ be the $n \times
n$ minor of $A$ got by removing the $i$th column.
$\Delta (t)$ divides $\det A_i$ so if $d = \max \{ \deg (\det A_i) \}$
then $\deg \Delta (t) \leq d$.  If each entry of $A$ has degree
$l$ or less then $\deg (\det A_i) \leq nl$, so all that remains
is to bound the degrees of the entries of $A$.

Let us consider an entry, $a_{ij} = \varphi \circ \psi (\frac{\partial
r_j}{\partial g_i} )$, of $A$.  As noted above, $r_j$ is a word in the
$g_k$'s of length at most $3(n+1)$.  Applying the rules above it is clear
that $\frac{\partial}{\partial g_i}r_j$ is a linear combination
of words in the $g_k$'s with lengths bounded by $3(n+1)$.  The map
$\varphi \circ \psi$ takes a word in the $g_k$'s to $t^a$ where
$a \in \mathbb{Z}$ is the number of times the path represented by
the word winds around the knot $L$.  Choose the integers $\nu_k$ so
that $\varphi \circ \psi (g_k ) = t^{\nu_k}$.  If $| \nu_k | \leq \nu$
for all $k$ then $\deg (a_{ij}) \leq \nu 3(n+1)$.

In the proof of Theorem \ref{hypthm} we saw that $|\nu_j| \leq 3^n$.
Hence,
\begin{equation*}
\deg (a_{ij}) \leq 3^n \cdot 3(n+1),
\end{equation*}
and hence
\begin{equation*}
\deg \Delta(t) \leq 3^{n+1}n(n+1).
\end{equation*}
\end{proof}

In Lemma \ref{decomposition_triangulation_lemma} we saw that a piece in the JSJ decomposition of
our manifold has at most $25 \cdot 4c = 100c$ tetrahedra.  The above
lemma tells us
that if the piece is a knot complement then its Alexander polynomial
has width at most $((100c)^2 + 100c)3^{100c+1}$.  The Alexander polynomial
of a $(u,v)$-torus knot has degree $(u-1)(v-1)$ (See
\cite[example 9.15]{Burde}).  Clearly $u,v \geq 2$,
whence, $(u-1) \leq (u-1)(v-1)$ and $(v-1) \leq (u-1)(v-1)$.
From these inequalities we get that $uv = (u-1)(v-1) + (u-1) + (v-1) + 1
\leq 3(u-1)(v-1) + 1$.
It follows that
\begin{equation}
uv \leq ((100c)^2 + 100c)3^{100c+2} + 1
\label{toruseq}
\end{equation}




\subsection{Covering Seifert fibered pieces}
\label{covering_seifert_fibered_pieces_subsection}

Now that we have bounded the orders of the cone points of in the
base orbifolds of our Seifert fibered pieces we may proceed to produce
the desired covers of these pieces.

\begin{Lem}
Let $N$ be Seifert fibered with base orbifold a disk with
two cone points of order $u$ and $v$.  Then for each prime $p > 3$,
$N$ has a cover, $\widetilde{N}$, with at most $2uvp^2$
sheets in which each boundary component of $\widetilde{N}$ is the
noncyclic cover with $p^2$ sheets of a boundary component of $N$.
\label{seiflem}
\end{Lem}

\begin{proof}
Let $N$ satisfy the hypotheses of the lemma, and 
let $F$ be the base orbifold of $N$.
$F$ is a disk with two cone points of order $u$ and $v$.

Glue a disk with one cone point of order $p> 3$ to $F$ to get a sphere 
$F'$ with 3 cone points with orders $u$, $v$, and $p$.
The orbifold, $F'$, is 
hyperbolic, and by \cite{Edmonds} has a finite orbifold cover 
which is a manifold. In fact, \cite{Edmonds} gives
such a cover, $\widetilde{F'}$,  with at most $2 \cdot
\mathrm{LCM} (u,v,p)$ sheets.
By removing open disk neighborhoods of each of 
the points of $\widetilde{F'}$ mapping to the cone point
of $F'$ with order $p$ 
we get a cover $\widetilde{F}$ of $F$.  By construction each 
boundary component of $\widetilde{F}$ is the $p$-fold cover of the 
boundary component of $F$.  This shows that $N$ has a cover, $N_0$, 
with at most $2uvp$ sheets whose base orbifold is the manifold,
$\widetilde{F}$,
and whose $S^1$ fibers map homeomorphically to the regular $S^1$ fibers of 
$N$.  We may then take $\widetilde{N}$ to be the $p$-fold cover of 
$N_0$ whose base space is again $\widetilde{F}$ and whose $S^1$ fibers 
are $p$-fold covers of the $S^1$ fibers of $N_0$.  Clearly 
$\widetilde{N}$ is the desired cover and has at most $2uvp^2$ 
sheets.
\end{proof}

Lemma \ref{seiflem}
and Inequality (\ref{toruseq}) combine to give the following theorem:

\begin{Thm} Let $N$ be a Seifert fibered
piece in the JSJ decomposition of the complement of a nontrivial knot,
$K$, with a diagram with $c$ crossings.  Let $p$ be any prime
greater than $3$.  Then 
$N$ has a cover $\widetilde{N}$ with
$((20000c^2 + 200c)3^{100c+2} + 2)p^2$ sheets or less whose
boundary components are the
noncyclic covers of order $p^2$ of the boundary components
of $N$. 
\label{seifthm}
\end{Thm}

Notice that the bound given in Theorem \ref{seifthm} is
exponential in the crossing number. It could be made polynomial
if it were known that the crossing number of a satellite
knot cannot be less than the crossing number of its companion.
It is conjectured that this should be true, but
it has remained unproven since Schubert introduced the
notion of satellite knots (See \cite[Problem 1.67]{Kirby}).




\section{Assembling the Covering Space}
\label{assembling_the_covering_space_section}

Now that we have produced the desired covers for the geometric
pieces of our knot complement, we must show that they can be
assembled to produce a cover of the entire complement.  This will be
done exactly as in \cite[section 2]{Hempel}.  Note that this
cover will in general not be regular.

\begin{figure}[ht]
$$
\setlength{\unitlength}{0.05in}
\put(-0.4,-1.5){\footnotesize \bf $\partial M$}
\put(43.4,2.45){\footnotesize \bf $M_0$}
\put(26.6,1.75){\footnotesize \bf $M_1$}
\put(9.45,1.19){\footnotesize \bf $M_2$}
\put(17.45,-2.05){\footnotesize \bf $T_2$}
\put(35.3,-1.28){\footnotesize \bf $T_1$}
\put(55,2.8){\footnotesize \bf $M$}
\put(45.85,21.7){\footnotesize \bf $\widetilde{M}_0$}
\put(25.5,21){\footnotesize \bf $\widetilde{M}_1$}
\put(6.5,20.3){\footnotesize \bf $\widetilde{M}_2$}
\put(17.1,15.6){\footnotesize \bf $\widetilde{T}_2$}
\put(34.4,15.7){\footnotesize \bf $\widetilde{T}_1$}
\put(14.8,21){\scriptsize \bf $\cong$}
\put(20,25){\begin{rotate}{55} \scriptsize  \bf $\cong$ \end{rotate}}
\put(19.48,18.21){\begin{rotate}{-45} \scriptsize  \bf $\cong$ \end{rotate}}
\put(32,25.4){\begin{rotate}{-45} \scriptsize  \bf $\cong$ \end{rotate}}
\put(32.7,18.5){\begin{rotate}{35} \scriptsize  \bf $\cong$ \end{rotate}}
\put(37.9,27){\begin{rotate}{50} \scriptsize  \bf $\cong$ \end{rotate}}
\put(37.4,21.7){\scriptsize \bf $\cong$}
\put(37,17.2){\begin{rotate}{-50} \scriptsize  \bf $\cong$ \end{rotate}}
\includegraphics[scale=0.7]{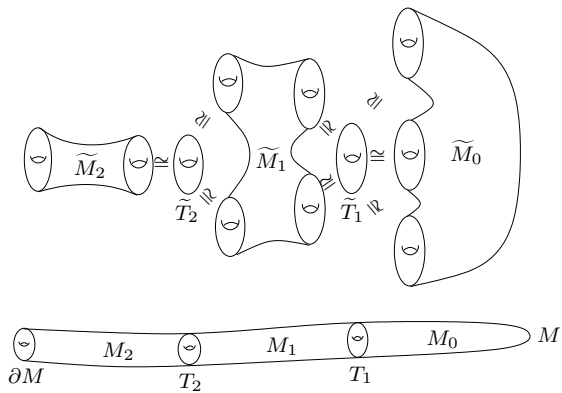}
$$
\caption{Pieces of the cover}
\label{cover_fig}
\end{figure}

\begin{figure}[ht]
$$
\setlength{\unitlength}{0.05in}
\put(5.4,63.5){\footnotesize \bf $\widetilde{M}^5_2$}
\put(5.4,48.45){\footnotesize \bf $\widetilde{M}^4_2$}
\put(5.4,41.45){\footnotesize \bf $\widetilde{M}^3_2$}
\put(5.4,26.4){\footnotesize \bf $\widetilde{M}^2_2$}
\put(5.4,17.5){\footnotesize \bf $\widetilde{M}^1_2$}
\put(5.4,2.45){\footnotesize \bf $\widetilde{M}^0_2$}
\put(14.6,56.5){\footnotesize \bf $\widetilde{M}^2_1$}
\put(14.6,34.1){\footnotesize \bf $\widetilde{M}^1_1$}
\put(14.8,11.0){\footnotesize \bf $\widetilde{M}^0_1$}
\put(24.8,50.2){\footnotesize \bf $\widetilde{M}^1_0$}
\put(25,16.6){\footnotesize \bf $\widetilde{M}^0_0$}
\includegraphics[scale=0.7]{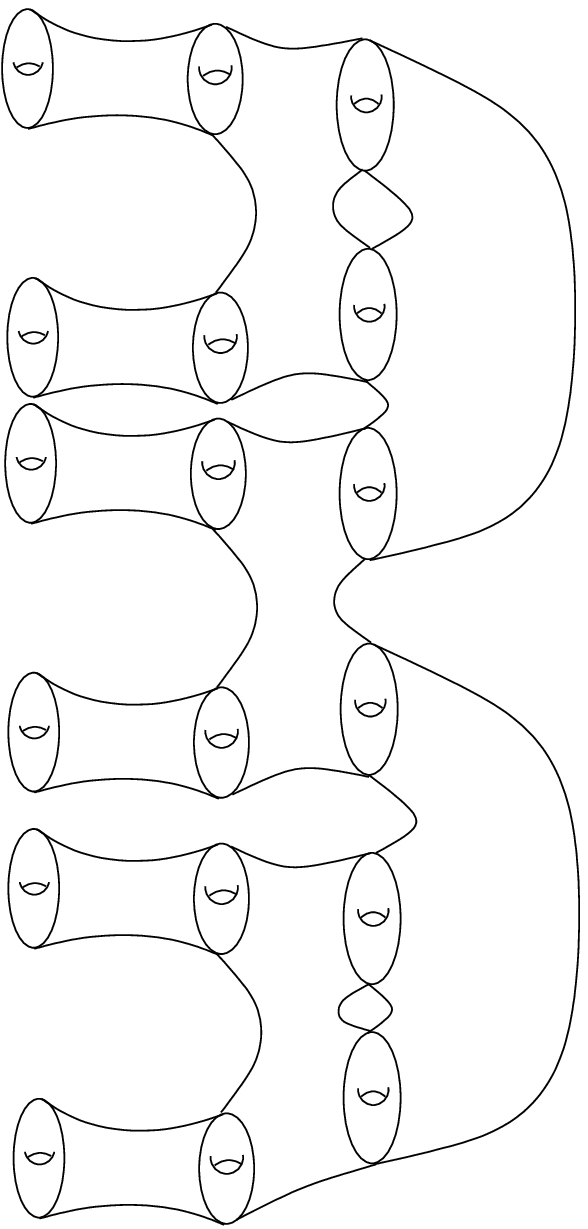}
$$
\caption{An assembled cover, $\widetilde{M}$, of $M$}
\label{assemble_fig}
\end{figure}

Recall that $K$ is a decompositionally linear knot in $S^3$, $M = S^3 - K$
its complement, and $\{ T_i \}_{i=1}^r$ a set of
tori cutting $M$ into geometric pieces $ M_0, M_1, \dots, M_r$
(See the bottom of Figure \ref{cover_fig}).
Fix a prime $p$. Each boundary component of $M_i$ is a torus with
fundamental group isomorphic to $\mathbb{Z} \times \mathbb{Z}$.
This group has a characteristic subgroup $P$ of index $p^2$
generated by $(p,0)$ and $(0,p)$.  For each $i$, let $\widetilde{T}_i$
be the cover of $T_i$ associated to the subgroup of $\pi_1(T_i)$
corresponding to $P \unlhd \mathbb{Z} \times \mathbb{Z}$.
Suppose for each $i$ we produce
a cover, $\widetilde{M}_i$, of $M_i$ such that the boundary
components of $\widetilde{M}_i$ are all covers equivalent to
$\widetilde{T}_j$ for some $j$ (see Figure \ref{cover_fig}). Then
by taking sufficiently many copies, $\widetilde{M}^k_i$, of
$\widetilde{M}_i$ we may assemble a cover $\widetilde{M}$ of
$M$ (see Figure \ref{assemble_fig}).
In fact if $m_i$ is the number of sheets in the
cover $\widetilde{M}_i$ then
we will need $\mathrm{LCM} ( \frac{m_0}{p^2},\frac{m_1}{p^2},
\dots, \frac{m_r}{p^2} ) / ( \frac{m_i}{p^2} )$
copies of $\widetilde{M}_i$, and $\widetilde{M}$ will
have $p^2 \cdot \mathrm{LCM}(\frac{m_0}{p^2},
\frac{m_1}{p^2}, \dots, \frac{m_r}{p^2})$ sheets.
It is immediate that $\widetilde{M}$ is not a cyclic
cover of $M$ since the boundary components of
$\widetilde{M}$ are not cyclic covers the boundary
of $M$.

We will now find a prime $p$ for which we can produce
such a set of coverings $\{ \widetilde{M}_i \}_{i=0}^r$.
Let us assume that our knot, $K$, has a diagram with $c$
crossings and hence by Lemma \ref{complement_triangulation_lemma}
its complement,
$M$, may be triangulated with $t \leq 4c$ tetrahedra.  Corollary \ref{deccor}
showed that $r \leq 2t \leq 8c$ and that the $M_i$'s can all be
ideally triangulated with a total of $25t \leq 100c$ or less tetrahedra.
If $M_i$ satisfies the hypotheses of Lemma \ref{homology_lemma}
then there is a number $w_i$ with $w_i \leq \frac{c}{2}$ such that
for every prime, $p$, not dividing $w_i$, Lemma \ref{homology_lemma}
gives a cover, $\widetilde{M}_i$, with $p^2$ sheets.
If $M_i$ does not satisfy the hypotheses of Lemma \ref{homology_lemma}
and is hyperbolic then Theorem \ref{hypthm} implies that there is a
number $B_i \in \mathbb{N}$ with $B_i \leq D(100c)$ such that for
every prime not dividing $B_i$ there is a cover,
$\widetilde{M}_i$, of $M_i$ which has
at most $p^{3\left( 2^{4n+4} \left(8n^2 + 4n \right)^{
(4n+4)2^{4n+4}}\right)}$ sheets
all of whose boundary tori are the noncyclic $p^2$ cover
of a boundary torus of $M_i$.  If $M_i$ is a torus knot complement
then Theorem \ref{seifthm} gives such a cover $\widetilde{M}_i$ with at
most $((20000c^2 + 200c)3^{100c+2} + 2)p^2$ sheets for all primes
$p \geq 3$. We wish to find a prime $p$ not dividing $B_i$ for
any $1 \leq i \leq r$.  The function $D$ is super exponential
so if $M_i$ has $t_i$ tetrahedra and $\sum t_i \leq 100c$ then
$\prod D(t_i)$ will be greatest when all tetrahedra are in
a single $M_i$.  Hence,
$$\prod w_i\prod B_i \leq  \left( {\textstyle \frac{c}{2}}
\right)^{8c}D(100c).$$
The product of all the primes less than $x > 2$ is at least
$e^{x/87}$ (See \cite[page 85]{Hua})
so there must be a prime $p$ less than
$87 \left(\log (D(100c)) + 8c \log  \frac{c}{2} \right)$
which does not divide any $B_i$ or $w_i$.

Each geometric piece $M_i$ has a cover
$\widetilde{M}_i$ whose boundary components are the noncyclic
$p^2$-sheeted covers of the boundary components of $M_i$.

From copies of these pieces we may assemble
a cover $\widetilde{M}$ of $M$ with at most
$p^2 \cdot \mathrm{LCM}(\frac{m_0}{p^2},
\frac{m_1}{p^2}, \dots, \frac{m_r}{p^2})$ sheets.
\begin{eqnarray}
\lefteqn{p^2 \cdot \mathrm{LCM} \left(\frac{m_0}{p^2},
\frac{m_1}{p^2}, \dots, \frac{m_r}{p^2} \right) } \nonumber \\
& \leq & m_0 m_1 \cdots m_r / p^{2r-2} \nonumber \\ [-3mm]
& & \\ [-3mm]
& \leq &  \left(87 \left(\log (D(100c)) + 8c \log {\textstyle
\frac{c}{2}} \right) \right)^{24c \left( 2^{4n+4} \left(8n^2 + 4n \right)^{
(4n+4)2^{4n+4}}\right)}. \nonumber
\end{eqnarray}
Thus we have completed the proof of Theorem \ref{mainthm}.

\section{Conclusions}

\label{concl_section}
We have shown that every decompositionally linear $c$-crossing knot has a finite non-cyclic
cover with at most $\Phi (c)$ sheets.  It should be noted that
this is very much a worst case result.  For example, the Alexander
polynomial of a knot determines when the fundamental group
of its complement surjects onto a dihedral group (See \cite{Burde} 14.8).

This result raises a number of questions. Firstly, it would
be nice to have such a bound for all knots.
The bound $\Phi (c)$ seems to be far from tight.  It would be
interesting to improve this bound.  From the other direction one
might try to produce lower limits for such a bound.  This could
be addressed by producing an infinite class of examples with a
large number of sheets in the smallest finite noncyclic cover
relative to the minimal crossing number.  From an algorithmic
point of view one might ask how to find and verify noncyclic
covers efficiently.

\appendix

\section{Appendices}

\subsection{An ideal triangulation of a knot complement}
\label{complement_triangulation_appendix}

For completeness we construct a triangulation of a knot complement
from a diagram of the knot.  This discussion is largely based on
the ideal triangulation algorithm in Jeffery Week's program, SnapPea.
However, we will present it from the viewpoint of standard spines.

\begin{proof}[Proof of Lemma \ref{complement_triangulation_lemma}]
Let $K$ be a knot in $S^3$ with a diagram with $c>0$ crossings.
We will produce a standard spine based on this projection
which will have less than $4c$ singular vertices.  This
spine will be dual to an ideal triangulation of the knot
complement with less than $4c$ ideal tetrahedra.
\begin{figure}[ht]
$$
\setlength{\unitlength}{0.05in}
\includegraphics[scale=1.0]{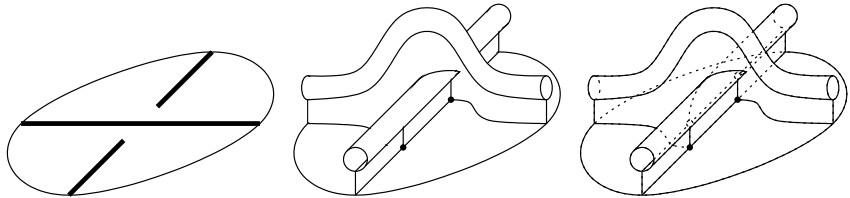}
$$
\caption{A knot projection and two views of
the corresponding spine $C'$}
\label{projspine_fig}
\end{figure}

We will assume that the projection of $K$ is a $4$-valent graph
in the equatorial 2-sphere, $S$, in $S^3$ with crossing
information at each vertex. Create
a spine, $C'$, as shown in Figure \ref{projspine_fig}.  Note
that $C'$ has exactly $4c$ singular
vertices.  For convenience we will assume that $K$ is
the core curve of the ``tubes'' in Figure \ref{projspine_fig}.
It should be clear that $C'$ satisfies the neighborhood condition.
As long as our projection contains at least one crossing,
the complement in $C'$ of its singular $1$-skeleton will be 
 a union of disks.  However,
$C'$ is not a spine of $M = S^3 - K$ since
two components of $M-C'$ are separated from $\partial
M$.

We now modify $C'$ to produce a standard spine
of $M$.  If one imagines $C'$ to be a soap-bubble film
the modifications we will make amount to ``popping''
two of the walls.
Let $C_1'$ be the singular 1-skeleton of $C'$.
Note that $C'$ divides $S^3$ into 3 components.  Two are
homeomorphic to open $3$-balls, and one is an open
regular neighborhood of $K$ which we
will denote by $N$.
Choose a disk, $D_1$, in $\partial N - C_1'$.
The closure of $D_1$ in $S^3$ contains a maximal open
arc, $\lambda$, such that $D_1 \cup \lambda$
is an open annulus.  Let $D_2$ be the
open disk in $C'-(C_1' \cup D_1)$ whose closure contains
$\lambda$.  Set $C'' = C' - (D_1 \cup
\lambda \cup D_2)$, and let $C_1''$ be its
singular 1-skeleton.  Choose an open disk
$D_3$ in $S - C_1''$.  Let $C = C'' - D_3$.

The spine, $C$, is, in fact, a standard spine for $M$. 
To prove this we first show that $M$ collapses to $C$.
It should be apparent than $\bar{N} - K$ collapses
to $\partial N$.  The removal of $(D_1 \cup
\lambda \cup D_2)$ from $C'$ joins $N-K$ to
a region homeomorphic to an open solid torus
along an open annulus.  The closure in $M$ of this new
region collapses
to it's boundary.  Finally, the removal of $D_3$
joins an open 3-ball to this region along an open
disk.  This joined region is $M-C$ and we
see that its closure in $M$ collapses to $C$.

Now we must show that $C$ is a standard spine.
The modification from $C'$ to $C$ preserves the
neighborhood condition.  Let
$C_1$ be the singular $1$-skeleton of $C$.
We will see that
$C''-C_1''$ is a disjoint union of open disks and then
that $C-C_1$ is as well.

First consider the effect
of removing $(D_1 \cup \lambda \cup D_2)$ from
$C'$ to get $C''$.  The equatorial 2-sphere, $S$, intersects
$C_1'$ in a graph which agrees with
the actual knot projection except near crossings.
In fact, $C_1'$ and the knot projection cut
$S$ into nearly identical 2-disks.  In $S$, the removal of
$(D_1 \cup \lambda \cup D_2)$ eliminates the
arc in $C_1'$ corresponding to the projection
of $\lambda$ onto $S$.  This has the effect
of joining certain disks in $S-C_1'$.  These
joined regions will all be disks because the
knot diagram will still be connected after
the removal of any one over-arc.
The other change after the removal of 
$(D_1 \cup \lambda \cup D_2)$ is that
the vertical walls at both ends of $\lambda$
will be joined to the pieces of the tunnel
on the other side of the wall from $\lambda$.
This amounts to gluing a disk to a disk along
an arc in their boundaries.  The results are
still disks.  Hence, $C''- C_1''$ is a disjoint
union of open disks.

Now consider the effect of removing $D_3$ from
$C''$ to get $C$.  Here each vertical wall surrounding
$D_3$ will be joined to the disk in $S$ on the 
other side of the wall from $D_3$.  These disks
are joined along single arcs in their boundary;
therefore, they glue together to form disks.
Consequently, $C-C_1$ is a union of open disks.

Let $C_0$ be the set of singular vertices of
$C$. I claim that $C_1 - C_0$ must be a disjoint union of
open arcs. If not then $C_1 - C_0$ contains an $S^1$. The
spine, $C$, is connected
and $C - C_1$ is composed
entirely of disks so each disk must, in fact, have this $S^1$ as its boundary.
This is impossible because the only such spine satisfying the neighborhood
condition is composed of three disks glued along their boundary.
This is not a spine
of the complement of a knot in the $3$-sphere or the complement
of a knot in a solid torus. Consequently, $C$ is a
standard spine of $M-K$.

The standard spine, $C$, has strictly
fewer singular vertices than $C'$, so $C$ has less than
$4c$ singular vertices.  It follows that there is a
dual ideal triangulation of $M$ with less than $4c$ ideal
tetrahedra.
\end{proof}

\subsection{Homology calculations}
\label{homology_calculations_appendix}

Here we will calculate the first homology group of a piece
in the JSJ decomposition of a knot complement. In section 
\ref{decomposition_of_a_knot_complement_subsection}
we saw that $M_0$ will always be
the complement of a knot in $S^3$.  A well-known Mayer-Vietoris
argument demonstrates that $H_1(M_0) \cong \mathbb{Z}$. The following
lemma gives the homology of pieces with two boundary components.
\begin{Lem}
Suppose $N$ is the complement of an open, regular neighborhood of
a knot $L$ in a solid torus.  Suppose further that $L$ has winding
number $w$ in the solid torus.  Then $H_1(N)$ has abelian 
presentation 
\begin{eqnarray*}
H_1(N) & = &  \Big\langle [m_0],[m_1],[l_0],[l_1] \Big| [l_1] =
w \cdot [m_0], [l_0] = w \cdot [m_1] \Big\rangle \\
 & = & \big\langle [m_0],[m_1] \big\rangle.
\end{eqnarray*}
\label{homology_calculation_lemma}
\end{Lem}
\begin{proof}
As in the statement of the lemma, let $N$ be the complement of
an open, regular neighborhood of a knot $L$ in a solid torus $V$.
Suppose that $L$ has winding number $w$ in the solid torus.
Fix an embedding of $V$ in $S^3$. Let $\partial_0 N$ and
$\partial_1 N$ be the boundary components
of $N$ as in section \ref{decomposition_of_a_knot_complement_subsection},
and let $m_k$ and $l_k$ be a meridian and longitude in $\partial_k N$.
For any curve $\alpha$ in $N$, let $[\alpha] \in H_1(N)$ be its homology
class.
Set $X$ to be an open regular neighborhood of $L$ in $V$ which intersects
$\mathring{N}$ in a regular neighborhood of  $\partial_0 N$.
Since $X \cup \mathring{N} = V$ we have the reduced homology
Mayer-Vietoris sequence
$$0 \to H_1(\partial_0 N) \to H_1(N) \oplus H_1(X) \to H_1(V) \to 0.$$
This exact sequence shows that  $H_1(N)$ is generated by
$H_1(\partial_0 N)$
and $H_1(V)$.  Of course $H_1(V)$ is generated by $H_1(\partial_1 N)$,
so we conclude that $H_1(N)$ is generated by $H_1(\partial_0 N)$
and $H_1(\partial_0 N)$.

The curve $l_1$ bounds a disk, $D$, in $V$.  One should observe
that $D \cap N$ demonstrates that
$[l_1] = w \cdot [m_0]$.  Longitude $l_0$ is homologically unlinked
with $L$ in $S^3$. Let $F$ be a Seifert surface for $l_0$ in $S^3$.
The surface, $F \cap N$, demonstrates that $[l_0] = w \cdot [m_1]$.
It follows that $\langle [m_0],[m_1] \rangle = H_1(N)$.
One sees that $[m_1]$ has infinite order in $H_1(N)$ by noting that
its image in $H_1(V)$ has infinite order.  Similarly, the image of
$[m_0]$ has infinite
order in $H_1(S^3-L)$ which implies that $[m_0]$ has infinite order
in $H_1(N)$.  Clearly, $[m_0]$ has trivial image in $H_1(V)$, and
$[m_1]$
has trivial image in $H_1(S^3-L)$.  This shows that $[m_1]$ and
$[m_0]$ are independent over the integers.
\end{proof}

\bibliographystyle{plain}
\bibliography{knotcov.bib}

\begin{thebibliography}{10}

\bibitem{Burde}
Gerhard Burde and Heiner Zieschang.
\newblock {\em Knots}.
\newblock de Gruyter Studies in Mathematics 5. Walter de Gruyter \& Co.,
  Berlin, 1985.

\bibitem{Casler}
B.~G. Casler.
\newblock An imbedding theorem for connected $3$-manifolds with boundary.
\newblock {\em Proc. Amer. Math. Soc.}, 16:559--566, 1965.

\bibitem{Culler}
M.~Culler and P.~B. Shalen.
\newblock Bounded, separating, incompressible surfaces in knot manifolds.
\newblock {\em Invent. Math.}, 75(3):537--545, 1984.

\bibitem{Dube}
Thomas~W. Dub\'{e}.
\newblock The structure of polynomial ideals and {G}r\"{o}bner bases.
\newblock {\em SIAM J. Comput.}, 19(4):750--775, 1990.

\bibitem{Edmonds}
Allan~L. Edmonds, John~H. Ewing, and Ravi~S. Kulkarni.
\newblock Torsion free subgroups of {F}uchsian groups and tessellations of
  surfaces.
\newblock {\em Invent. Math.}, 69(3):331--346, 1982.

\bibitem{Hass}
Joel Hass and Jeffrey~C. Lagarias.
\newblock The number of reidemeister moves needed for unknotting.
\newblock {\em J. Amer. Math. Soc.}, 14(2):399--428, 2001.

\bibitem{Hempel}
John Hempel.
\newblock Residual finiteness for $3$-manifolds.
\newblock In {\em Combinatorial group theory and topology (Alta, Utah, 1984)},
  Ann. of Math. Stud., 111, pages 379--396. Princeton Univ. Press, Princeton,
  NJ, 1987.

\bibitem{Hua}
Loo~Keng Hua.
\newblock {\em Introduction to number theory}.
\newblock Springer-Verlag, Berlin-New York, 1982.
\newblock Translated from the Chinese by Peter Shiu.

\bibitem{Kawauchi}
Akio Kawauchi.
\newblock {\em A survey of knot theory}.
\newblock Birkh\"{a}user Verlag, Basel, 1996.
\newblock Translated and revised from the 1990 Japanese original by the author.

\bibitem{Kirby}
Rob Kirby.
\newblock Problems in low-dimensional topology.
\newblock Published electronically at {\tt
  \verb|http://math.berkeley.edu/~kirby/problems.ps.gz|}, 1995.

\bibitem{Mignotte}
Maurice Mignotte and Doru {\c{S}}tef\u{a}nescu.
\newblock {\em Polynomials. An algorithmic approach}.
\newblock Springer Series in Discrete Mathematics and Theoretical Computer
  Science. Springer-Verlag, Singapore, 1999.

\bibitem{Murasugi}
Kunio Murasugi.
\newblock {\em Knot theory and its applications}.
\newblock Springer Series in Discrete Mathematics and Theoretical Computer
  Science. Birkh\"{a}user Boston, Inc., Boston, MA, 1996.
\newblock Translated from the 1993 Japanese original by Bohdan Kurpita.

\bibitem{NS}
Walter~D. Neumann and Gadde~A. Swarup.
\newblock Canonical decompositions of $3$-manifolds.
\newblock {\em Geom. Topol.}, 1:21--40, 1997.

\bibitem{Petronio}
Carlo Petronio.
\newblock Ideal triangulations of link complements and hyperbolicity equations.
\newblock {\em Geom. Dedicata}, 66(1):27--50, 1997.

\bibitem{Piergallini}
Riccardo Piergallini.
\newblock Standard moves for standard polyhedra and spines.
\newblock {\em Rend. Circ. Mat. Palermo (2) Suppl.}, 18:391--414, 1988.

\bibitem{Rolfsen}
Dale Rolfsen.
\newblock {\em Knots and links}.
\newblock Mathematics Lecture Series, 7. Publish or Perish, Inc., Houston, TX,
  1990.
\newblock Corrected reprint of the 1976 original.

\bibitem{Schmidt}
Wolfgang~M. Schmidt.
\newblock Heights of algebraic points.
\newblock In {\em Number theory and its applications (Ankara, 1996)}, Lecture
  Notes in Pure and Appl. Math., 204, pages 185--225. Dekker, New York, 1999.

\bibitem{Schubert}
Horst Schubert.
\newblock \"{U}ber eine numerische knoteninvariante.
\newblock {\em Math. Z.}, 61:245--288, 1954.
\newblock (German).

\bibitem{TV}
V.~G. Turaev and O.~Ya Viro.
\newblock State sum invariants of $3$-manifolds and quantum $6j$-symbols.
\newblock {\em Topology}, 31(4):865--902, 1992.

\bibitem{Zhang}
Shou-Wu Zhang.
\newblock Personal communication.
\newblock email: {\tt szhang@math.columbia.edu}.

\end{thebibliography}

\end{document}